# UNIFORM CONVERGENCE OF EXACT LARGE DEVIATIONS FOR RENEWAL REWARD PROCESSES[1]

By Zhiyi Chi

*University of Connecticut*

Let $(X_n, Y_n)$ be i.i.d. random vectors. Let $W(x)$ be the partial sum of $Y_n$ just before that of $X_n$ exceeds $x > 0$. Motivated by stochastic models for neural activity, uniform convergence of the form $\sup_{c \in I} |a(c,x) \Pr\{W(x) \geq cx\} - 1| = o(1)$, $x \to \infty$, is established for probabilities of large deviations, with $a(c,x)$ a deterministic function and $I$ an open interval. To obtain this uniform exact large deviations principle (LDP), we first establish the exponentially fast uniform convergence of a family of renewal measures and then apply it to appropriately tilted distributions of $X_n$ and the moment generating function of $W(x)$. The uniform exact LDP is obtained for cases where $X_n$ has a subcomponent with a smooth density and $Y_n$ is not a linear transform of $X_n$. An extension is also made to the partial sum at the first exceedance time.

**1. Introduction.**

1.1. *Background.* Let $(X_n, Y_n)$ be i.i.d. $\sim (X, Y) \in \mathbb{R}^2$. For $x \geq 0$, let

$$N(x) = \begin{cases} \max\left\{n : \sum_{i=1}^{k} X_i \leq x, \forall k \leq n\right\}, & \text{if } X_1 \leq x, \\ 0, & \text{otherwise}, \end{cases}$$

$$S(x) = \sum_{i=1}^{N(x)} X_i, \qquad W(x) = \sum_{i=1}^{N(x)} Y_i,$$

where $\sum_{i=1}^{n} x_i := 0$ if $n = 0$. It will hereafter be assumed that

(1.1) $\quad X$ and $Y$ are nondegenerate and $\Pr\left\{\sup_{n \geq 1} \sum_{i=1}^{n} X_i = \infty\right\} = 1.$

Received January 2005; revised December 2005.
[1]Supported in part by NIH Grants MH68028, DC007206.
*AMS 2000 subject classifications.* Primary 60F10; secondary 60G51.
*Key words and phrases.* Large deviations, renewal reward process, point process, continuous-time random walk.

The process $\{W(x), x \geq 0\}$ is often called a *continuous-time random walk* or *renewal reward process* (cf. [1] and references therein). There have been many studies of the LDPs for $N(x)$ and $W(x)$ and related issues for exit points, often in the more general context of Markov additive or renewal processes [7, 9, 10, 11, 15, 16, 18, 19, 21, 22, 23, 26].

In this article, we consider the uniform exact LDP for the random sum $W(x)$, which refers to asymptotics of the form

$$\sup_{c \in I} |a(c,x) \Pr\{W(x) \geq cx\} - 1| = o(1) \qquad \text{as } x \to \infty,$$

with $I$ being an open interval and $a(c,x)$ a deterministic function. Our interest in this problem originates from research in neuroscience. In electrophysiological studies, neuronal activity is recorded as sequences of pulses in voltage generated by neurons and can be modeled as point processes. In some cases, the neuronal activity can be regarded as consisting of consecutive "epochs" such that the activities therein are driven by the same underlying mechanism [4]. Each epoch yields a measurement, such as the number of pulses during its duration, and the cumulative measurement within a period of time is of interest. Let $X_k$ and $Y_k$ be the duration and measurement of the $k$th epoch, respectively. Then the overall measurement up to time $x$ is $W(x)$. It is useful to understand the probabilistic properties of $W(x)$, such as the fluctuations in the number of pulses generated by a neuron over the long term. In the simplest models, the epochs are independent of each other, which gives rise to the problem at hand. Although $X_n$ is positive for neural activity, in general this is not a required assumption.

The exact LDP is known for the case where $c$ is fixed (cf. [15, 26]) and for other related cases, such as the degenerate case $Y \equiv 1$ [8]. However, the available results are not sufficient in some applications. For instance, when dealing with the number of pulses in neural activity, $X_n$ may vary continuously while $Y_n$ are integers. Since $W(x) \geq cx$ is equivalent to $W(x) \geq \lceil cx \rceil = c^*x$, with $c^*$ varying according to $x$, it is necessary to consider the convergence for $c^*$ in a neighborhood of $c$ instead of for $c$ alone.

1.2. *Overview.* Our "plan of attack" is as follows. We choose the uniform exact LDP of [3] as the basic tool, which requires some careful analysis of the moment generating function of $W(x)$. Henceforth, denote

(1.2) $\quad M_x(z) = E[e^{zW(x)}], \qquad g_x(z) = g(z,x) = E[e^{zY}|X=x], \qquad z \in \mathbb{C}.$

It is seen that for $t \in \mathbb{R}$, $M_x(t) = \Pr\{X > x\} + \int_{-\infty}^{x} g_u(t) M_{x-u}(t) F(du)$. We have to make sure that asymptotics of the form $M_x(t) \sim a(x,t) e^{b(t)x}$ hold uniformly well for $t \in I$, where $a(x,t)$ and $b(t)$ are deterministic functions and $I \neq \varnothing$ is an open interval. Suppose that $E[e^{tY - h(t)X}] \equiv 1$ for some function $h(t)$. Letting $\phi_x(t) = \frac{M_x(t)}{e^{h(t)x}}$ then yields renewal equations $\phi_x(t) = \psi_x(t) + \int_{-\infty}^{x} \phi_{x-u}(t) F_t(du)$, where $F_t(du) = \frac{g_u(t)}{e^{h(t)u}} F(du)$ is a proper probability measure. By the standard theory [2, 5], for each $t$, $\phi_x(t)$ converges and so $M_x(t) \sim \phi_x(t) e^{h(t)x}$. However, it is necessary to show uniform convergence for $\phi_x(t)$. One approach is to first establish the uniform convergence of the renewal measures associated with $F_t$. Indeed, if $h(t)$ is smooth, then the family of $F_t$ is smooth in a certain sense, which facilitates the establishment of the uniform convergence of the renewal measures. Next, provided $\psi_x(t)$ are uniformly integrable, the uniform convergence of $\phi_x(t)$ can be obtained. We can then check necessary conditions on $M_x$ to obtain the desired result for $W(x)$. In particular, we need to establish uniform bounds on $\frac{M_x(t+si)}{M_x(t)}$. Suffice it to say, our derivation of the bounds critically depends on the uniform asymptotics of $M_x(t)$, the smoothness of $F_t$ and the assumption that $Y$ is *not* a deterministic affine function of $X$.

*Notation.* Recall that a *connected component* of $A \subset \mathbb{R}$ is a maximal nonempty interval in $A$ and the *interior* of $A$ is $A^o = \{x : (x - r, x+) \subset A \text{ for some } r > 0\}$. The Lebesgue measure will be denoted by $\ell$. Define $a \vee b = \max\{a, b\}$, $a \wedge b = \min\{a, b\}$. By convention, we let $\inf \varnothing = \infty$ and $\sup \varnothing = -\infty$.

1.3. *Main results.* To carry out the plan outlined above, we shall first generalize the convergence of renewal measures associated with one measure [17, 24, 25] to exponentially fast uniform convergence for a family of measures. This can be done in a setting not specific to the uniform exact LDP. For a probability measure $p$ on $\mathbb{R}$, denote by $\mu_p$ its mean and by $N_p = \sum_{n=0}^{\infty} p^{n*}$ the associated renewal measure, where $p^{n*}$ is the $n$-fold convolution of $p$, with $p^{0*} := \delta(x)$. Denote by $\text{sppt}(p)$ the support of $p$. If $X \sim p$, also let $\text{sppt}(X) = \text{sppt}(p)$. The next result generalizes those in [12], which are restricted to $X \geq 0$, and gives no bounds on the convergence of densities. For uniform convergence with power rates, see [13, 14].

THEOREM 1.1 (Uniform convergence of renewal measures). *Let $\mathcal{M}$ be a family of probability measures on $\mathbb{R}$ such that*

(1.3) $$\inf_{p \in \mathcal{M}} \mu_p > 0,$$

(1.4) $$M_\tau := \sup_{p \in \mathcal{M}} E_p[e^{\tau|X|}] < \infty \quad \text{for some } \tau > 0.$$

*Suppose that each $p \in \mathcal{M}$ is a mixture of two probability measures $\Phi_p$ and $\Psi_p$, where*

(1.5) $$p = (1 - \lambda_p)\Phi_p + \lambda_p \Psi_p, \qquad 0 \leq \lambda_p < 1,$$

such that $\Phi_p$ has a density $\phi_p$ and

(1.6) $$\gamma := \sup_{p \in \mathcal{M}} \lambda_p < 1.$$

Furthermore, suppose that there exists $T \in (0, \infty)$ such that

(1.7) $$\phi_p \in C_0^2((-T, T)) \qquad \forall p \in \mathcal{M}$$

(1.8) $$\sup_{p \in \mathcal{M}} \int_{-T}^{T} |\phi_p''(x)| \, dx < \infty.$$

Then $N_p = Q_p + \bar{N}_p$ such that $Q_p$ has a density $q_p$, $\sup_{p \in \mathcal{M}} \bar{N}_p(\mathbb{R}) < \infty$ and for all $0 < r \ll 1$,

(1.9) $$A_1 := \sup_{x > 0, p \in \mathcal{M}} \{e^{rx} |q_p(x) - \mu_p^{-1}|\} + \sup_{x > 0, p \in \mathcal{M}} \{e^{rx} q_p(-x)\} < \infty,$$

(1.10) $$A_2 := \sup_{x > 0, p \in \mathcal{M}} \{e^{rx} \bar{N}_p((-\infty, -x] \cup [x, \infty))\} < \infty.$$

We will need two corollaries of Theorem 1.1. Corollary 1.2 is based on Corollary 1.1 and is the one which will be used directly in obtaining the uniform exact LDP.

COROLLARY 1.1. *Let $z_p(x)$ be functions on $\mathbb{R}$. Suppose that there exists $\eta > 0$ such that $L := \sup_{x \in \mathbb{R}, p \in \mathcal{M}} \{e^{\eta |x|} |z_p(x)|\} < \infty$. Then for each $p$, there exits a unique solution $Z_p$ to*

(1.11) $$Z_p(x) = z_p(x) + \int_{-\infty}^{\infty} Z_p(x - y) p(dy)$$

*such that* $\sup_{x \in \mathbb{R}} \{e^{-\varepsilon x} Z_p(x)\} < \infty, \forall \varepsilon > 0.$

*Let $u_p = \mu_p^{-1} \int z_p(x) \, dx$. Let $r, A_1, A_2$ be the same as in Theorem 1.1 and let $\alpha = \frac{1}{2} \min(\eta, r)$. Then there exists a continuous function $C(\cdot)$ such that*

(1.12) $$\sup_{x > 0, p \in \mathcal{M}} (e^{\alpha x} |Z_p(x) - u_p|) + \sup_{x < 0, p \in \mathcal{M}} (e^{\alpha |x|} |Z_p(x)|) \leq C(A_1, A_2, L, \eta, r).$$

COROLLARY 1.2. *Let $z_p(x)$ be defined on $x \geq 0$ and $\eta \in (0, \tau)$ such that $L_1 := \sup_{x \geq 0, p \in \mathcal{M}} \{e^{\eta x} |z_p(x)|\} < \infty$, where $\tau$ is as in (1.4). Then*

(1.13) $$Z_p(x) = z_p(x) + \int_{-\infty}^{x} Z_p(x - y) p(dy) \qquad \forall x \geq 0, \; p \in \mathcal{M},$$

has at most one solution satisfying $L_2(\varepsilon) := \sup_{x\geq 0, p \in \mathcal{M}} \{e^{-\varepsilon x} Z_p(x)\} < \infty$, $\forall \varepsilon > 0$. If the solution exists, then as $x \to \infty$,

$$(1.14) \quad Z_p(x) \to u_p := \frac{1}{\mu_p}\left[\int_0^\infty z_p(x)\,dx - \int_{-\infty}^0 \left(\int_0^{|y|} Z_p(x)\,dx\right)p(dy)\right]$$

and for $\alpha = \frac{1}{2}\min(\eta, r)$, $L = L_1 \vee (L_2(\eta)M_\tau)$ and $C(\cdot)$ as in Corollary 1.1,

$$(1.15) \quad \sup_{x>0, p\in\mathcal{M}} (e^{\alpha x}|Z_p(x) - u_p|) \leq C(A_1, A_2, L, \eta, r).$$

Conditions (1.6) and (1.7) in Theorem 1.1 imply that Stone's decomposition can apply uniformly well for $p \in \mathcal{M}$ [25]. Conditions (1.7) and (1.8) imply that $\phi_p$ have bounded oscillations, so they cannot cluster around a lattice distribution. The conditions also imply that the characteristic functions of $\phi_p$ are uniformly integrable and can hence be inverted, which facilitates the analysis of the tail of $\phi_p$. As indicated earlier, the condition (1.4) can be met by tilting a measure with a smoothly parameterized family of functions.

Returning to the LDP for $W(x)$, define

$$\mathcal{D}_t = \{s \in \mathbb{R} : E[e^{tY-sX}] \leq 1\} \qquad \forall t \in \mathbb{R},$$
$$h(t) = \inf \mathcal{D}_t, \qquad h^*(\nu) = \sup_t\{\nu t - h(t)\},$$
$$\beta_X = \limsup_{x\to\infty}(1/x)\log \Pr\{X > x\} \leq 0.$$

For what values of $c$ can a nontrivial LDP be expected for $\Pr W(x) \geq cx$? If $X > 0$ and $Y$ are independent, then the LLN implies that $c > c_0 := \frac{EY}{EX}$, due to $N(x) \approx \frac{x}{EX}$ and $W(x) \approx N(x)EY \approx c_0 x$. In general, the answer depends on the properties of $h(t)$. In the following theorems, $\tau_0 > 0$ is a basic assumption in order for a nontrivial LDP to arise for $c = h'(t)$ when $t$ takes values across a neighborhood of $\tau_0$. Note that by the assumptions of the theorems, $Y$ is not a deterministic linear function of $X$. Therefore, the case $Y \equiv 1$ studied in [8] is not covered here.

THEOREM 1.2 (Uniform exact LDP, nonlattice case). *Let $\tau_0 > 0$ and*

$$(1.16) \qquad h(\tau_0) \in (\beta_X, \infty), \qquad \inf_{s\in\mathbb{R}} E[e^{\tau_0 Y + sX}] < 1.$$

*Suppose the following conditions (1)–(4) are satisfied.*

(1) *The law of $X$ can be decomposed into the sum of two nonnegative measures $\Phi$ and $\Psi$ such that $\Phi(\mathbb{R}) > 0$ with a density $\phi \in C^2(\mathbb{R})$.*

(2) *There exists $k \in \mathbb{N}$ and $\eta_0 > 0$ such that*

$$(1.17) \qquad E|E[e^{is(Y_1+\cdots+Y_k)}|X_1+\cdots+X_k]| < 1 \qquad \forall s \in \mathbb{R}\setminus\{0\},$$

$$(1.18) \qquad E[e^{-h(\tau_0)X+\tau_0 Y+\eta_0(|X|+|Y|)}] < \infty.$$

(3) *With the same $\eta_0$ as in (2),*

$$g(t,x), \frac{\partial g(t,x)}{\partial x}, \frac{\partial^2 g(t,x)}{\partial x^2} \in C((\tau_0 - \eta, \tau_0 + \eta) \times \mathrm{sppt}(\Phi)^o).$$

(4) *In the case where $\Pr\{X \geq 0\} < 1$,*

(1.19) $$E[e^{q(\tau_0 Y - (0 \wedge h(\tau_0))X)} | X > 0] < \infty \qquad \forall q > 0.$$

*Then there exists $I = (\tau_0 - \varepsilon, \tau_0 + \varepsilon) \neq \varnothing$ such that $h \in C^\infty(I)$, $E[e^{tY - h(t)X}] \equiv 1$ and $h''(t) > 0$ for $t \in I$ and*

(1.20) $$\sup_{t \in I} \left| \frac{t}{\phi(t)} \sqrt{2\pi x h''(t)} e^{xh^*(h'(t))} \Pr\{W(x) \geq xh'(t)\} - 1 \right| = o(1)$$

*as $x \to \infty$, where*

(1.21) $$0 < \phi(t) := \frac{B(X,t)}{E[X e^{tY - h(t)X}]} < \infty,$$

*with $0 < B(X,t) < \infty$ given by*

$$B(X,t) = \int_0^\infty \frac{\Pr\{X > x\}}{e^{h(t)x}} \, dx - \int_{-\infty}^0 \int_0^{|x|} \frac{g_x(t) M_u(t)}{e^{h(t)(x+u)}} \, du F(dx).$$

*In particular, if $X \geq 0$ a.s., then*

$$B(X,t) = \begin{cases} h(t)^{-1}(1 - E[e^{-h(t)X}]), & \text{if } h(t) \neq 0, \\ E[X], & \text{if } h(t) = 0. \end{cases}$$

COROLLARY 1.3. *Under the assumptions of Theorem 1.2, given $a, b \in \mathbb{R}$,*

(1.22) $$\Pr\{W(x+a) \geq xh'(\tau_0) + b\} \sim \frac{\phi(\tau_0) e^{-xh^*(h'(\tau_0)) + ah(\tau_0) - b\tau_0}}{\tau_0 \sqrt{2\pi x h''(\tau_0)}}.$$

If $Y$ has a finite moment generating function, then the results generalize to the sums at the first passage time without much extra effort. Let

$$T(x) = \min\left\{n : \sum_{i=1}^n X_i > x\right\}, \qquad \bar{W}(x) = \sum_{i=1}^{T(x)} Y_i.$$

COROLLARY 1.4. *Suppose that (1.16) and conditions (1)–(3) of Theorem 1.2 hold and that $E[e^{tY}] < \infty$ for all $t$. Additionally, if $\Pr\{X \geq 0\} < 1$, then instead of (1.19), suppose that*

(1.23) $$E[e^{-q\min(0,h(\tau_0))X} \mathbf{1}\{X > 0\}] < \infty \qquad \forall q > 0.$$

*Then with $W(x)$ being replaced by $\bar{W}(x)$ and $\phi(t)$ replaced by*

$$0 < \bar{\phi}(t) = \frac{\bar{B}(X,t)}{E[X e^{tY - h(t)X}]} < \infty,$$

(1.22) still holds. The function $\bar{B}(X,t)$ is defined as

$$\bar{B}(X,t) = \int_0^\infty \int_u^\infty \frac{g_x(t)}{e^{h(t)(x+u)}} F(dx)\,du - \int_{-\infty}^0 \int_0^{|x|} \frac{g_x(t)E[e^{t\bar{W}(u)}]}{e^{h(t)(x+u)}}\,du F(dx).$$

In particular, if $X \geq 0$ a.s., then

$$\bar{B}(X,t) = \begin{cases} h(t)^{-1}E[e^{tY}(1-e^{-h(t)X})], & \text{if } h(t) \neq 0, \\ E[Xe^{tY}], & \text{if } h(t) = 0. \end{cases}$$

Now consider the case where $Y$ is lattice-valued. Recall that a random variable $\xi$ is said to have *span* $d > 0$ if $d = \max\{t > 0 : \Pr\{\xi \in t\mathbb{Z}\} = 1\}$.

THEOREM 1.3 (Uniform exact LDP, lattice case). *Let $Y$ be lattice-valued with span $d$. Suppose that all of the conditions of Theorem 1.2 are satisfied except for (1.17). Instead, assume that there exists $k \in \mathbb{N}$ such that*

(1.24) $\quad E|E[e^{is(Y_1+\cdots+Y_k)}|X_1+\cdots+X_k]| < 1 \quad \forall s \in (0, \pi/d].$

*Then with $\tau_0 > 0$ and $\phi$ as in Theorem 1.2 and $\{a\} := \lceil a \rceil - a$, as $x \to \infty$,*

(1.25) $\Pr\{W(x) \geq xh'(\tau_0)\} \sim \dfrac{\phi(\tau_0)d}{\sqrt{2\pi x h''(\tau_0)}} \dfrac{e^{-xh^*(h'(\tau_0)) - d\tau_0\{xh'(\tau_0)/d\}}}{1 - e^{-\tau_0 d}}.$

To see the relevance of the condition $\tau_0 > 0$, consider again the case where $X$ and $Y$ are independent, with $EX > 0$. Then by $E[e^{tY}]E[e^{-h(t)X}] = 1$, $h'(0) = \frac{EY}{EX}$. By the strict convexity of $h$, $h'(\tau_0) > \frac{EY}{EX}$. As mentioned just before Theorem 1.2, this gives rise to a nontrivial LDP.

Some comments on (1.17) are in order. One of its implications is that for any $a, b \in \mathbb{R}$ and $d > 0$, $\Pr\{Y \in aX + b + d\mathbb{Z}\} < 1$. Apparently, if $E|E[e^{isY}|X]| < 1$, $\forall s \neq 0$, then (1.17) is satisfied with $k = 1$. Another condition that implies (1.17) is as follows.

PROPOSITION 1.1. *Suppose that $X$ has a continuous density and that there exists a piecewise continuous function $f$ on $\mathrm{sppt}(X)$ such that $Y - f(X)$ is constant or lattice valued. Suppose that the graph of $f$ is not a set of parallel straight line segments, that is, for any $c \in \mathbb{R}$, $f(x) - cx$ is not piecewise constant. Then $(X, Y)$ satisfies condition (1.17).*

Finally, in order to attain the uniform exact LDP, we need $\phi(t) > 0$. If $X \geq 0$, then it is easy to see that $\phi(t) > 0$. The case where $\Pr\{X \geq 0\} < 1$ is more involved. The condition (1.19) is imposed to ensure that in order for $\phi(\tau_0) > 0$.

In what follows, Section 2 gives examples of application of the uniform exact LDP of $W(x)$. In Section 3, Theorem 1.1 and its corollaries are proved. In Section 4, the main theorems concerning the uniform exact LDP are proved. Section 5 proves corollaries and related results for the uniform exact LDP. Section 6 collects auxiliary technical details for the main results.

## 2. Examples.

EXAMPLE 2.1. Consider the example on neural activity in Section 1. Assume that the duration $X$ of an epoch has density $e^{-x}\mathbf{1}\{x>0\}$ and that within each epoch, the neural activity has the same dynamics. Specifically, let $f \geq 0$ be a piecewise $C^1$ function on $(0,\infty)$ with $\int f > 0$. Given an epoch of the form $[\eta, \eta + X]$, the point process $\xi$ therein is such that $\xi - \eta$ is a Poisson process with density $f(t)\mathbf{1}\{t \in [0, X]\}$. Let $Y$ be the number of points in the epoch. Then $Y|X \sim \text{Poisson}(F(X))$, with $F(x) = \int_0^x f$. Suppose that $\lim_{x\to\infty} \frac{F(x)}{x} = \nu \geq 0$. Let $\psi(z) = \int_0^\infty e^{z(F(x)-\nu x)}\,dx$. For any $x > 0$, $W(x)$ is lattice-valued with span $d = 1$. Define $a = -\infty$ if $\psi(z) > 1\ \forall z \in \mathbb{R}$, otherwise $a = \sup\{z \in \mathbb{R} : \psi(z) \leq 1\}$.

We show the exact LDP (1.25) holds for $\tau_0 > 0$ with $e^{\tau_0} - 1 > a$. We need to check conditions (1) and (3) of Theorem 1.2, (1.16), (1.18), and (1.24). Condition (1) is clear. Given $X = x > 0$, $E[e^{tY}|X = x] = \exp\{F(x)(e^t - 1)\}$. Let $\{x_n\}$ be the set of points where $f$ fails to be $C^1$. Let $\Phi$ be the density of $X$ restricted to the exterior of a neighborhood of $\{x_n\}$. Then condition (3) is satisfied. Now $\forall s \in (0, \pi]$, $E|E[e^{isY}|X]| = E|e^{F(X)(e^{is}-1)}| = E[e^{F(X)(\cos s - 1)}] < 1$. Therefore, (1.24) holds.

For any $s$, $E[e^{\tau_0 Y - sX}] = b(s) := \int_0^\infty e^{F(x)(e^{\tau_0}-1)-(s+1)x}\,dx$. Observe that $b(s)$ is continuous and strictly decreasing on $\{b(s) < \infty\}$. Since $\frac{F(x)}{x} \to \nu$ as $x \to \infty$, for $s > s_0 := \nu(e^{\tau_0}-1) - 1 \geq 1$, $b(s) < \infty$, and $b(s) \to 0$ as $s \to \infty$. On the other hand, as $s \downarrow s_0$, $b(s) \uparrow \psi(e^{\tau_0}-1)$ and by $e^{\tau_0}-1 > a$, $\lim_{s\downarrow s_0} b(s) > 1$. As a result, there exists a unique $s = h(\tau_0) > s_0$ with $E[e^{\tau_0 Y - sX}] = 1$. It is then easy to see that there exists $\eta > 0$ such that $E[e^{(\tau_0+\eta)Y-sX+\eta X}] < \infty$. Both (1.16) and (1.18) therefore hold.

By (1.21), $\phi(t) = [(1+h(t))\int_0^\infty x e^{F(x)(e^t-1)-(h(t)+1)x}\,dx]^{-1}$ for $t > 0$ with $e^t - 1 > a$. Now suppose that $f(x) \equiv 1$. In this case, $a = -\infty$ and $W(x)$ is the number of points in $[0, (x-X) \vee 0]$ from a Poisson process $V(x)$ with density 1. We obtain $h(t) = e^t - 1$, $\phi(t) = e^{-t}$ and $h^*(\nu) = \nu \ln \nu - \nu + 1$ for $t > 0$, $\nu > 1$. For any $T > 0$, letting $\tau_0 = \ln T > 0$ in (1.25) yields

$$\Pr\{W(x) \geq Tx\} \sim \frac{e^{-x(T\ln T - T + 1) - \ln T\{Tx\}}}{\sqrt{2\pi Tx}(T-1)} = \frac{T^{-\lceil Tx\rceil}e^{(T-1)x}}{\sqrt{2\pi Tx}(T-1)}.$$

It is worth comparing the exact LDP of $W(x)$ with that of $V(x)$. For $x \gg 0$, the "cut-off" interval $[(x-X) \vee 0, x]$ is a very small fraction of $[0, x]$. Does the random cut-off have any effect? By Theorem 3.5 of [3], we have

$$\Pr\{V(x) \geq Tx\} = \Pr\{N(x) \geq \lceil Tx \rceil\} \sim \frac{e^{-x(\Lambda_x^*)'(\lceil Tx\rceil/x)}}{\sqrt{2\pi x \Lambda_x''(\tau_x)}} \frac{1}{1 - e^{-\tau_x}}, \qquad T > 1,$$

where $\Lambda_x(t) = \frac{1}{x}\log E[e^{tN(x)}] = e^t - 1$, $\Lambda_x^*(\nu) = \nu \ln \nu - \nu + 1$ for $\nu > 1$ and $\tau_x$ is the unique solution to $\Lambda_x'(\tau_x) = \frac{\lceil Tx\rceil}{x}$. It follows that the effect can be quantified as $\Pr\{V(x) \geq Tx\} \sim T\Pr\{W(x) \geq Tx\}$.

Now let $Y = \mathbf{1}\{X > M\}$, where $M > 0$ is a constant. Then $W(x)$ is the number of epochs in $[0, x]$ with duration longer than $M$. Conditioning on $X_1 + X_2 = u \in (0, 2M)$, $\zeta = Y_1 + Y_2$ is binary with $\Pr\{\zeta = 0\} \in (0, 1)$. Therefore, $|E[e^{is(Y_1+Y_2)}|X_1 + X_2 = u]| < 1$, $\forall s \in (0, \pi]$ and condition (1.24) is satisfied. For $t > 0$, $s = h(t)$ is the solution to $se^{Ms} = e^{-M}(e^t - 1)$. It follows that for any $c \in (e^{-M}, \frac{1}{M})$, there exists a unique $\tau_0 > 0$ with $h'(\tau_0) = c$. The exact LDP for $\Pr\{W(x) \geq cx\}$ can then be obtained.

EXAMPLE 2.2. Fix $a > 0$. Let $X \sim N(0,1)$ and $Y = \mathbf{1}\{X \geq a\} - \mathbf{1}\{X < a\}$. Regard $S_n = X_1 + \cdots + X_n$ as a random walk. Given $x > 0$, $N(x) < \infty$ a.s. and $W(x)$ is the difference between the number of steps with size greater than $a$ and the number with size less than $a$ before the random walk crosses $x$ for the first time. Note that $EN(x) = \infty$ and hence that $W(x)$ is not integrable. Let $\Phi(x) = \Pr\{X \leq x\}$. For $t > 0$ and $s \in \mathbb{R}$,

$$E[e^{tY-sX}] = e^{s^2/2} \underbrace{\left[e^t(1 - \Phi(s+a)) + e^{-t}\Phi(s+a)\right]}_{A(t)}.$$

Then $h(t) = -\sqrt{2 \ln \frac{1}{A(t)}} < 0$ is well defined only for $t \in (0, \ln \frac{\Phi(a)}{1-\Phi(a)})$. Now $h(t) > \beta_X = \lim_{x \to \infty} \frac{1}{x} \Pr\{X \geq x\} = -\infty$ and condition (1.19) is satisfied. As in Example 2.1, $W(x)$ satisfies the exact LDP in (1.25) with $d = 1$.

EXAMPLE 2.3. Suppose that $X$ and $Y$ are independent, each has a $C^2$ density such that $EX \geq EY = 1$, $\operatorname{ess\,inf} X = \operatorname{ess\,inf} Y = 0$ and

$$t_0 = \sup\{t : E[e^{tY}] < \infty\} > 0, \qquad \lim_{t \to t_0-} E[Ye^{tY}] = \infty,$$

where $\operatorname{ess\,inf} X := \sup\{x : \Pr\{X > x\} = 1\}$. Let $\xi_0 = (0,0)$ and for $n \geq 1$, let $\xi_n = (\xi_{n,1}, \xi_{n,2}) := (\sum_{i=1}^n X_i, \sum_{i=1}^k Y_i)$.

Given $x > 0$, let $I_x = [0, x] \times [0, x]$. Consider the probability that the first $\xi_n$ outside $I_x$ is a certain distance from the upper-right corner $(x, x)$. Let $T(x) = \min\{n : \xi_n \notin I_x^o\}$. Given $\lambda \in (0, 1)$ and $M \in \mathbb{R}$, let $A_{\lambda, x, M} = \{\xi_{T(x),1} \leq \lambda x + M\}$. Observe that, for $0 < y < x$, $\xi_{T(x),1} \leq y \iff$ "for some $k \geq 1$, $\sum_{i=1}^k X_i \leq y$ and $\sum_{i=1}^k Y_i > x$" $\iff W(y) > x$.

As $X \geq 0$ is nondegenerate, for $t > 0$, there exists a unique $h(t)$ with $E[e^{tY - h(t)X}] = E[e^{tY}]E[e^{-h(t)X}] = 1$. Then $h'(t) = \frac{a_t}{b_t}$, where $a_t = \frac{E[Ye^{tY}]}{E[e^{tY}]}$ is strictly increasing and $b_t = \frac{E[Xe^{-h(t)X}]}{E[e^{-h(t)X}]}$ strictly decreasing. If $h(t_0-) = \infty$, then as $t \uparrow t_0$, $b_t \to \operatorname{ess\,inf} X = 0$ and hence $h'(t) \geq \frac{a_0}{b_t} \to \infty$. On the other hand, if $h(t_0-)$ is finite, then as $t \uparrow t_0$, $E[e^{tY}] = \frac{1}{E[e^{-h(t)X}]}$ is bounded and hence, by $E[Ye^{tY}] \to \infty$, $h'(t) \geq \frac{a_t}{b_0} \to \infty$. In either case, $h'(t) \to \infty$. Since $h'(0) = \frac{EY}{EX} \leq 1$, there is a unique $\tau_\lambda$ with $h'(\tau_\lambda) = \frac{1}{\lambda}$. So by Corollary 1.3,

$$\Pr\{A_{\lambda,x,M}\} = \Pr\{W(\lambda x + M) \geq x\}$$
$$= \Pr\{W(\lambda x + M) \geq \lambda x h'(\tau_\lambda)\} \sim \frac{\phi(\tau_\lambda) e^{Mh(\tau_\lambda)}}{\tau_\lambda \sqrt{2\pi \lambda x h''(\tau_\lambda)}} e^{-\lambda x h^*(\lambda^{-1})},$$

where $\phi(\tau_\lambda) = \frac{1 - E[e^{-h(\tau_\lambda)X}]}{h(\tau_\lambda) E[Xe^{-h(\tau_\lambda)X}] E[e^{tY}]}$.

**3. Exponentially fast uniform convergence of renewal measures.** This section proves Theorem 1.1. For any $\sigma$-finite measure $p$, let $\hat{p}$ denote its characteristic function.

PROOF OF THEOREM 1.1. Since $\hat{N}_p = 1 + \hat{p}\hat{N}_p$ and $\hat{p} = (1 - \lambda_p)\hat{\Phi}_p + \lambda_p \hat{\Psi}_p$, following [25], $\hat{N}_p = (1 - \lambda_p \hat{\Psi}_p)^{-1}[1 + (1 - \lambda_p)\hat{\Phi}_p \hat{N}_p]$ and

$$(3.1) \qquad N_p = (1 - \lambda_p)\Phi_p * \bar{N}_p * N_p + \sum_{n=0}^{\infty} \lambda_p^n \Psi_p^{n*} = Q_p + \bar{N}_p.$$

Because $\Phi_p$ has a density $\phi_p \in C_0^2((-T, T))$ with $T > 0$ independent of $p$ and because $1 - \lambda_p > 0$, $Q_p$ has a density $q_p \in C$. Also, note that $\bar{N}_p(\mathbb{R}) = \frac{1}{1-\lambda_p}$.

We first show that (1.10) holds for $0 < r \ll 1$. By $\lambda_p \Psi_p \leq p$ and $\gamma = \sup_{p \in \mathcal{M}} \lambda_p < 1$, for any $\varepsilon > 0$, $x > 0$ and $p \in \mathcal{M}$, we have

$$\bar{N}_p((-\infty, -x] \cup [x, \infty)) \leq S_p + \sum_{n \geq \varepsilon x} \lambda_p^n \leq S_p + \frac{\gamma^{\varepsilon x}}{1 - \gamma}$$

where $S_p = \sum_{n < \varepsilon x} \Pr\{|\sum_{i=1}^n X_i| \geq x\}$ with $X_i$ i.i.d. $\sim p$. By (1.4), there exists $\tau > 0$, such that $M_\tau = \sup_{p \in \mathcal{M}} E_p[e^{\tau|X|}] < \infty$. Then by Chernoff's inequality,

$$S_p \leq \sum_{n \leq \varepsilon x} E_p[e^{\tau(nX - x)} + e^{-\tau(nX + x)}] \leq 2\varepsilon x e^{x(\varepsilon \log M_\tau - \tau)}.$$

Thus, given $\varepsilon \in (0, \frac{\tau}{\log M_\tau})$, (1.10) holds for $0 < r < \min(\varepsilon \log \frac{1}{\gamma}, \tau - \varepsilon \log M_\tau)$.

The rest of the proof is devoted to (1.9). We need a suitable spectral representation of $q_p(x)$. By $\bar{N}_p(\mathbb{R}) = \frac{1}{1-\lambda_p}$, $\chi_p = (1 - \lambda_p)\Phi_p * \bar{N}_p$ is a probability measure. Then, as in [25],

$$(3.2) \qquad q_p(x) = \frac{1}{2\mu_p} + \frac{1}{2\pi} \int_{-\infty}^{\infty} \Re\left(e^{-ix\theta} \frac{\hat{\chi}_p(\theta)}{1 - \hat{p}(\theta)}\right) d\theta,$$

$$(3.3) \qquad \frac{1}{2\pi} \int_{-\infty}^{\infty} \Re\left(e^{-ix\theta} \frac{\hat{\chi}_p(\theta)}{-i\theta}\right) d\theta = \frac{1}{2} - \chi_p([x, \infty)).$$

For completeness, proofs for (3.2) and (3.3) are given in the Appendix.

Because $\text{sppt}(\Phi_p) \subset [-T, T]$, the tail probability of $\chi_p$ is the same as that of $\bar{N}_p$ up to a shift bounded by $T$ and a multiplicative factor bounded by $1 - \gamma$. Then by (1.10),

$$(3.4) \qquad \sup_{x > 0, p \in \mathcal{M}} \{e^{rx} \chi_p((-\infty, -x] \cup [x, \infty))\} < \infty \qquad \forall 0 < r \ll 1.$$

By condition (1.3), (3.2)–(3.4) and the self-conjugacy of $\hat{\chi}_p(\theta)$ and $\hat{p}(\theta)$,

$$(3.5) \qquad q_p(x) = \begin{cases} \dfrac{1}{\mu_p} + \dfrac{1}{2\pi} I_p(x) + R_p(x), & \text{if } x > 0, \\ \dfrac{1}{2\pi} I_p(x) + R_p(x), & \text{if } x < 0, \end{cases}$$

with $\sup_{x > 0, p \in \mathcal{M}} |e^{\delta x} R_p(x)| < \infty$, where

$$(3.6) \quad I_p(x) = \int_{-\infty}^{\infty} e^{-ix\theta} \hat{\chi}_p(\theta) K_p(\theta) \, d\theta, \qquad K_p(z) := \frac{1}{1 - \hat{p}(z)} - \frac{1}{-i\mu_p z}.$$

To continue, we need the next lemma, which will be proven in Section 6.1.

LEMMA 3.1. *There exists $\eta > 0$ such that for all $p \in \mathcal{M}$, $\hat{p}(z) := E_p[e^{izX}]$ is analytic on $D_\eta = \{z : |\Im(z)| \leq \eta\}$ and $\hat{p}(z) \neq 1$ for $z \neq 0$. Furthermore, $L_\eta := \sup_{z \in D_\eta, p \in \mathcal{M}} |K_p(z)| < \infty$.*

Let $A_p(s) = \lambda_p \int e^{sx} \Psi_p(dx)$. By (3.1), $|\hat{\bar{N}}_p(\theta - is)| \leq \sum_{n=0}^{\infty} A_p(s)^n$, $s \in \mathbb{R}$. Given $s \in (-\tau, \tau)$, since $0 < e^{|sx|} - 1 \leq |sx| e^{|sx|} \leq \frac{|s| e^{\tau |x|}}{\tau - |s|}$,

$$0 \leq \lambda_p \int e^{|sx|} \Psi_p(dx) - \lambda_p \leq \int (e^{|sx|} - 1) p(dx) \leq \frac{|s| M_\tau}{\tau - |s|},$$

yielding $0 \leq A_p(s) \leq \lambda_p + \frac{|s| M_\tau}{\tau - |s|}$. Thus one can choose $0 < \eta \ll 1$ as in Lemma 3.1 such that $\sup_{p \in \mathcal{M}, |s| \leq \eta} A_p(s) < 1$. As a result, $\hat{\bar{N}}_p(z)$ are analytic and uniformly bounded in $D_\eta$. By means of some computation,

$$(3.7) \qquad |\hat{\Phi}_p(\theta - is)| \leq e^{|s|T} \min(1, C\theta^{-2}) \qquad \forall \theta, s \in \mathbb{R}, p \in \mathcal{M},$$

where $C = \sup_{p \in \mathcal{M}} \int |\phi_p''|$. It follows that $\hat{\chi}_p = (1 - \lambda_p) \hat{\Phi}_p \hat{\bar{N}}_p$ is analytic in $D_\eta$ and $(1 + \theta^2) \sup_{p \in \mathcal{M}, |s| \leq \eta} |\hat{\chi}_p(\theta - is)| < \infty$. Thus there exists a constant $C_1 > 0$ such that $\int_{-\infty}^{\infty} |\hat{\chi}_p(\theta - is)| \, d\theta < C_1$ for $s \in [-\eta, \eta]$ and $p \in \mathcal{M}$.

Let $r \in (0, \eta)$. By Lemma 3.1, $K_p(z)$ is analytic on $D_r$. For $x > 0$, apply Lemma 3.1 and Cauchy's contour integral to $e^{-izx} \hat{\chi}_p(z) K_p(z)$ along the path $\Im(z) = 0$ and $\Im(z) = -r$. Then by (3.6),

$$I_p(x) = \int_{-\infty}^{\infty} e^{-x(r + i\theta)} \hat{\chi}_p(\theta - ir) K_p(\theta - ir) \, d\theta$$

and hence $|I_p(x)| \leq L_r e^{-xr} \int_{-\infty}^{\infty} |\hat{\chi}_p(\theta - ir)| \, d\theta \leq C_1 L_r e^{-xr}$. Similarly, for $x < 0$, $|I_p(x)| \leq C_1 L_r e^{xr}$. So by (3.5), (1.9) holds. □

PROOF OF COROLLARY 1.1. For $p \in \mathcal{M}$, $Z_p(x) = \int_{-\infty}^{\infty} z_p(x-y) N_p(dy)$ is a solution to (1.11). For $x > 0$,

$$Z_p(x) - u_p = \int_{-\infty}^{\infty} z_p(x-y) \bar{N}_p(dy) + \int_{-\infty}^{\infty} z_p(x-y) q_p(y) \, dy$$
$$- \frac{1}{\mu_p} \int_{-\infty}^{\infty} z_p(x-y) \, dy$$

and hence

$$|Z_p(x) - u_p| \leq \int_{x/2}^{\infty} |z_p(x-y)| \bar{N}_p(dy) + \int_{-\infty}^{x/2} |z_p(x-y)| \bar{N}_p(dy)$$
$$+ \int_{x/2}^{\infty} |z_p(x-y)(q_p(y) - \mu_p^{-1})| \, dy$$
$$+ \int_{-\infty}^{x/2} |z_p(x-y)| (q_p(y) + \mu_p^{-1}) \, dy.$$

The integrals on the right-hand side are bounded by $LA_2 e^{-rx/2}$, $LA_2 e^{-\eta x/2}$, $LA_1 e^{-rx/2}$ and $\frac{1}{\eta} L(A_1 + \frac{2}{\mu_p}) e^{-\eta x/2}$, respectively. A similar bound for $x < 0$ can be obtained. Therefore, $Z_p$ satisfies (1.12).

It remains to show that the solution is unique for each $p$. If $\text{sppt}(p) \in [0, \infty)$, this follows from [5], Lemma 4.1.I. In general, let $D(x)$ be the difference between two such solutions. Then $D(x) = \int_{-\infty}^{\infty} D(x-y) p^{n*}(dy)$. Given $\varepsilon > 0$, there exists $C > 0$ such that $D(x-y) \leq C e^{\varepsilon(x-y)}$ and hence $|D(x)| \leq \int_{-\infty}^{\infty} C e^{\varepsilon(x-y)} p^{n*}(dy) = C e^{\varepsilon x} (E_p[e^{-\varepsilon X}])^n$. By (1.3) and (1.4), $E_p[e^{-\varepsilon X}] < 1$ if $\varepsilon \ll 1$. Letting $n \to \infty$ yields $D(x) = 0$. The uniqueness is thus proved. □

PROOF OF COROLLARY 1.2. The uniqueness of the solution to (1.13) can be shown by following the proof of Corollary 1.1. Let $Z_p(x)$, $x \geq 0$, be the solution. Extend $Z_p(x)$ and $z_p(x)$ to $\mathbb{R}$ so that for $x < 0$, $Z_p(x) = 0$ and $z_p(x) = -\int_{-\infty}^{x} Z_p(x-y) p(dy)$. Then the renewal equation

$$Z_p(x) = z_p(x) + \int_{-\infty}^{\infty} Z_p(x-y) p(dy)$$

holds. By (1.4), for $x \leq 0$, $|z_p(x)| \leq \int_{-\infty}^{x} L_2(\eta) e^{\eta(x-y)} p(dy) \leq L_2(\eta) e^{\eta x} M_\tau$. Thus $\sup_{x \in \mathbb{R}, p \in \mathcal{M}} \{e^{\eta |x|} |z_p(x)|\} \leq L$. By Corollary 1.1, (1.15) holds. □

**4. Uniform exact LDP.** This section proves Theorems 1.2 and 1.3. We start with some basic results.

4.1. *Preparations.* By condition (1.18),

$$\varphi(t,s) := E[e^{tY+sX}] \in C^\infty((t,s): |t-\tau_0| < \eta_0, |s+h(\tau_0)| < \eta_0).$$

Then $f(s) := \varphi(\tau_0, s) \in C^\infty$, $f(-h(\tau_0)) = 1$. By condition (1.16), $\inf f < 1$. As $h(\tau_0) = -\sup\{s : f(s) \le 1\}$ and $f$ is convex, $f'(-h(\tau_0)) > 0$. Therefore,

$$\frac{\partial \varphi(\tau_0, -h(\tau_0))}{\partial s} = E[X e^{\tau_0 Y - h(\tau_0) X}] > 0.$$

Fix $\tau_0$ and $\eta_0$ as in Theorem 1.2. By the implicit function theorem (cf. [20]), there exist $0 < \varepsilon_0 < \tau_0 \wedge \eta_0$ and $0 < \eta < (h(\tau_0) - \beta_X) \wedge \eta_0$ such that $h$ is a smooth mapping from $I_0 = [\tau_0 - \varepsilon_0, \tau_0 + \varepsilon_0]$ into $[h(\tau_0) - \eta, h(\tau_0) + \eta]$,

(4.1) $$E[e^{tY - h(t)X}] \equiv 1$$

and

(4.2) $$E[X e^{tY+sX}] = \frac{\partial \varphi(t,s)}{\partial s} > 0 \qquad \text{on } I_0 \times [-h(\tau_0) - \eta, -h(\tau_0) + \eta].$$

Differentiate $E[e^{tY-h(t)X}] \equiv 1$ twice to obtain $h''(t) = \frac{E[(Y-h'(t)X)^2 e^{tY-h(t)X}]}{E[X e^{tY-h(t)X}]}$. By condition (1.17), $h''(t) > 0$ and hence $h(t)$ is strictly convex on $I_0$.

By (4.1), for any $t \in I_0$, define the following probability measure:

$$P_t(dx, dy) = e^{ty - h(t)x} P(dx, dy).$$

Denote by $E_t$ the expectation under $P_t$ and by $E$ that under $P$. Then $E_t[X] = E[X e^{tY-h(t)X}]$ and $\sup_{t \in I_0} E_t[X] < \infty$. Under $P_t$, the marginal of $X$ and the conditional measure of $Y$ given $X$ are, respectively,

(4.3)
$$F_t(dx) = g_x(t) e^{-h(t)x} F(dx),$$
$$P_t(dy|x) := P_t(dy|X=x) = \frac{e^{ty}}{g_x(t)} P(dy|x).$$

Note that for any function $f$, $E_t[f(X)] = E[e^{tY-h(t)X} f(X)]$. Define

(4.4) $$\psi_x(t) = \frac{\Pr\{X > x\}}{e^{h(t)x}}, \qquad \phi_x(t) = \frac{M_x(t)}{e^{h(t)x}}, \qquad k_x(t) = \frac{g_x(t)}{e^{h(t)x}}.$$

Then $F_t(dx) = k_x(t) F(dx)$. To apply Corollary 1.2 to the proof of Theorem 1.2, the following lemmas are needed.

LEMMA 4.1. *The family of probability measures $\{F_t(dx), t \in I_0\}$ satisfies conditions (1.3)–(1.8).*

LEMMA 4.2. *For all $x \geq 0$ and $t \in I_0$, $M_x(t) < \infty$. Define*

$$(4.5) \qquad \Lambda_x(t) = \frac{1}{x} \log M_x(t), \qquad \Lambda_x^*(\lambda) = \sup_{t \in \mathbb{R}}\{\lambda t - \Lambda_x(t)\}$$

*and define $\phi(t)$ by (1.21). Then there exists $I = [\tau_0 - \varepsilon, \tau_0 + \varepsilon] \subset I_0^o$ with $\varepsilon \in (0, \varepsilon_0 \wedge \tau_0)$ such that $\inf_{t \in I} h''(t) > 0$,*

$$(4.6) \qquad \limsup_{x \to \infty} \sup_{t \in I}\{x|\Lambda_x^{(n)}(t) - h^{(n)}(t)|\} < \infty, \qquad n = 0, 1, 2,$$

*and $\phi(t)$ is continuous on $I$, $\inf_{t \in I} \phi(t) > 0$, and there exists $\alpha > 0$ such that*

$$(4.7) \qquad \limsup_{x \to \infty} e^{\alpha x} \sup_{t \in I} |\phi_x(t) - \phi(t)| < \infty.$$

In the following proofs, the uniform exact LDP in [3] is the fundamental tool. However, it turns out that some of the conditions used in the main result of [3] are hard to verify for $W(x)$. We shall instead check the more basic conditions provided in that work.

4.2. *Nonlattice case.* PROOF OF THEOREM 1.2. Let $I$ be as in Lemma 4.2. To prove (1.20), it suffices to show that for any $x_n \to \infty$ and $\tau_n \in (\tau_0 - \varepsilon/2, \tau_0 + \varepsilon/2)$,

$$(4.8) \quad \Pr\{W(x_n) \geq x_n h'(\tau_n)\} \sim \frac{\phi(\tau_n)}{\tau_n \sqrt{2\pi x_n h''(\tau_n)}} e^{-x_n h^*(h'(\tau_n))}, \qquad n \to \infty.$$

For brevity, define $M_n(t) = M_{x_n}(t)$ and $\phi_n(t) = \phi_{x_n}(t)$. Let $\Lambda_n(t) = \frac{1}{x_n} \times \log M_n(t)$ and $\nu_n = h'(\tau_n)$. Note that $\alpha := \inf_{I_0} h''(t) > 0$. By (4.6),

$$C := \frac{1}{2\alpha} \limsup_n \left[x_n \sup_{t \in I} |\Lambda_n'(t) - h'(t)|\right] < \infty.$$

Let $\varepsilon_n = \frac{C}{x_n}$. For $n \gg 1$, $\tau_n \pm \varepsilon_n \in I$ and hence $\Lambda_n'(\tau_n + \varepsilon_n) - \nu_n \geq \Lambda_n'(\tau_n + \varepsilon_n) - h'(\tau_n + \varepsilon_n) + \alpha \varepsilon_n > 0$. Likewise, $\nu_n - \Lambda_n'(\tau_n - \varepsilon_n) > 0$. Since $\Lambda_n'$ is strictly increasing, there exists a unique $t_n^* \in I$ with

$$(4.9) \qquad \Lambda_n'(t_n^*) = \nu_n, \qquad |t_n^* - \tau_n| \leq C/x_n.$$

Define random variables $A_n^*$ and $U_n$ such that

$$\Pr\{A_n^* \in du - x_n \nu_n\} = \frac{\exp\{t_n^* u\}}{M_n(t_n^*)} \Pr\{W(x_n) \in du\},$$

$$U_n = \frac{A_n^*}{\sqrt{x_n \Lambda_n''(t_n^*)}}.$$

We will later show that there exist $\delta > 0$ and $N > 0$, such that

$$(4.10) \qquad f^*(s) := \sup_{n \geq N} \{|E[e^{isU_n}]|\mathbf{1}\{|s| \leq \delta\sqrt{x_n \Lambda_n''(t_n^*)}\}\} \in L^1,$$

$$(4.11) \quad \sup_{\eta < |s| < \lambda} \left|\frac{E[e^{(t_n^* + is)W(x_n)}]}{M_n(t_n^*)}\right| = o\left(\frac{1}{\sqrt{x_n}}\right) \qquad \forall \eta, \lambda > 0.$$

Equation (4.10) implies (2.7) in [3]. Meanwhile, it can be shown that $U_n \xrightarrow{d} N(0,1)$, which is the first conclusion of Lemma 3.1 of [3] (cf. [6], Theorem 3.7.4). Thus the first claim in the proof of Lemma 3.2 of [3] holds. As $0 < \liminf t_n^* < \limsup t_n^* < \infty$, $t_n^*$ satisfies (3.3) in [3] and (4.11) implies condition (c) of Lemma 3.2 of [3]. Thus that lemma holds for $W(x_n)$. By Theorem 3.3 of [3],

$$\Pr\{W(x_n) \geq x_n \nu_n\} \sim \frac{e^{-x_n \Lambda_n^*(\nu_n)}}{t_n^* \sqrt{2\pi x_n \Lambda_n''(t_n^*)}}, \qquad n \to \infty.$$

By (4.6), $\frac{\Lambda_n''(t_n^*)}{h''(\tau_n)} \to 1$. Note that $e^{-x_n \Lambda_n^*(\nu_n)} = e^{-x_n(\nu_n t_n^* - h(t_n^*))}\phi_n(t_n^*)$. Since $\phi \in C(I)$, by (4.7) and (4.9), $\lim \phi_n(t_n^*) = \phi(\tau_0)$. By Taylor expansion,

$$\nu_n t_n^* - h(t_n^*) = \nu_n \tau_n - h(\tau_n) - \tfrac{1}{2}h''(\xi_n)(t_n^* - \tau_n)^2$$
$$= h^*(\nu_n) - \tfrac{1}{2}h''(\xi_n)(t_n^* - \tau_n)^2, \qquad \text{some } \xi_n \in [\tau_n, t_n^*].$$

By (4.9), $\nu_n t_n^* - h(t_n^*) = h^*(\nu_n) + o(x_n^{-1})$ and hence (1.20) is proved. □

4.3. *Lattice case.* The proof for this case is based on the next lemma.

LEMMA 4.3. *Under the assumptions of Theorem* 1.3, *for all $x \gg 1$, $W(x)$ is lattice-valued with span $d$.*

PROOF OF THEOREM 1.3. As in the proof of Theorem 1.2, the bound in (4.10) still holds. On the other hand, following an argument similar to that in the proof for (4.11), it is not hard to see that for any $\eta \in (0, \pi/d]$,

$$\sup_{\eta < |s| \leq \pi/d} \left|\frac{E[e^{(t_x^* + is)W(x)}]}{M_x(t_x^*)}\right| = o\left(\frac{1}{\sqrt{x}}\right).$$

The rest of the proof follows that of Theorem 1.2. Let $\nu = h'(\tau_0)$ and $\nu_x = \frac{d}{x}\lceil \frac{x\nu}{d}\rceil$. As $|\nu_x - \nu| \leq \frac{d}{x}$, $\forall x \gg 0$, there exits a unique $t_x^*$ with $h'(t_x^*) = \nu_x$. By Theorem 3.5 of [3],

$$\Pr\{W(x) \geq x\nu\} = \Pr\left\{W(x) \geq d\left\lceil\frac{x\nu}{d}\right\rceil\right\}$$
$$= \Pr\{W(x) \geq x\nu_x\}$$
$$\sim \frac{d}{\sqrt{2\pi x \Lambda_x''(t_x^*)}}\frac{e^{-x\Lambda_x^*(\nu_x)}}{1 - e^{-t_x^* d}}.$$

By $|h'(t_x^*) - h'(\tau_0)| = |\nu_x - \nu| \leq \frac{d}{x}$, there exists $C > 0$ with $|t_x^* - \tau_0| \leq \frac{C}{x}$. Then $\Lambda_x''(t_x^*) \to h''(\tau_0)$, $e^{-x\Lambda_x^*(\nu_x)} = e^{-x(\nu_x t_x^* - \Lambda_x(t_x^*))} = e^{-x(\nu_x t_x^* - h(t_x^*))}\phi_x(t_x^*) \sim e^{-x(\nu_x t_x^* - h(t_x^*))}\phi(\tau_0)$ and

$$\nu_x t_x^* - h(t_x^*) = (\nu_x - \nu)(t_x^* - \tau_0) + (\nu_x - \nu)\tau_0 + \nu t_x^* - h(t_x^*)$$
$$= (\nu_x - \nu)\tau_0 + \nu\tau_0 - h(\tau_0) + O(x^{-2}).$$

Finally, note that $x(\nu_x - \nu)\tau_0 = d\tau_0\{x\nu/d\}$. This then completes the proof. □

4.4. *Uniform bounds for the tilted complex moment generating functions.* The rest of the section is devoted to (4.10) and (4.11). Let $\sigma_n = \sqrt{\Lambda_n''(t_n^*)}$. Then

$$E[e^{zU_n}] = \exp\left\{-\frac{\nu_n\sqrt{x_n}z}{\sigma_n}\right\}\frac{1}{M_n(t_n^*)}M_n\left(t_n^* + \frac{z}{\sqrt{x_n}\sigma_n}\right), \qquad z \in \mathbb{C}.$$

PROOF OF (4.10). The formula will follow if there exist $a, \delta > 0$ such that

$$\limsup_n \sup_{s\in\mathbb{R}}\{|E[e^{isU_n}]|\mathbf{1}\{|s| \leq \delta\sqrt{x_n}\sigma_n\}e^{as^2}\} < \infty.$$

Replacing $s$ with $\sqrt{x_n}\sigma_n s$, the supremum becomes

$$\sup_{|s|\leq\delta}\left\{\frac{|M_n(t_n^* + is)|}{M_n(t_n^*)}e^{ax_n\sigma_n^2 s^2}\right\}.$$

As $t_n^* \in I$ and $\Lambda_n''(t_n^*)$ is bounded away from 0 and $\infty$, (4.10) will follow if there exist $a > 0$ and $\delta > 0$ such that

$$\limsup_{x\to\infty}\sup_{t\in I, |s|\leq\delta}\left\{\frac{|M_x(t+is)|}{M_x(t)}e^{axs^2}\right\} < \infty.$$

Note that $M_x(t) = \phi_x(t)e^{h(t)x}$. Let $G_x(t+is) = e^{-h(t)x}E[e^{(t+is)W(x)}]$. By Lemma 4.2, $\sup_I|\frac{\phi_x(t)}{\phi(t)} - 1| \to 0$ and $\inf_I \phi(t) > 0$. Thus it suffices to show that

(4.12) $$\limsup_{x\to\infty}\sup_{t\in I, |s|\leq\delta}\{|G_x(t+is)|e^{axs^2}\} < \infty.$$

Define $w_k = \sum_{i=1}^k u_i$. By (4.3) and (4.4), for $s, t \in \mathbb{R}$,

$$G_x(t+is) = \psi_x(t) + \int_{-\infty}^x E_t[e^{isY}|u]G_{x-u}(t+is)F_t(du)$$
$$= \sum_{n=0}^\infty \underbrace{\int \psi_{x-w_n}(t)\mathbf{1}\{w_k \leq x, \forall k \leq n\}\prod_{j=1}^n E_t[e^{isY}|u_j]F_t(du_j)}_{J_n(s,t,x)}.$$

Note that $|J_n(s,t,x)| \le J_n(t,x) := J_n(0,t,x)$ and $\sum_{n\ge 0} J_n(t,x) = G_x(t) = \phi_x(t)$. First, we show that there exist $\eta > 0$ and $M > 0$ such that

$$(4.13) \quad \limsup_{x\to\infty} e^{\eta x}\left\{\sup_{t\in I}\sum_{n\le \eta x}|J_n(t,x)| + \sup_{t\in I}\sum_{n\ge Mx}|J_n(t,x)|\right\} < \infty.$$

By (1.18) and $\inf_I h(t) > \beta_X$, there exists $a \in (0, \frac{1}{2}(\inf_I h(t) - \beta_X))$ such that

$$L := \sup_{t\in I}\log E_t[e^{2a|X|}] < \infty,$$

$$K(\psi) := \sup_{t\in I, x\ge 0}|e^{2ax}\psi_x(t)| < \infty.$$

Then

$$|J_n(t,x)| \le K(\psi)\int e^{-2a(x-w_n)}\prod_{j=1}^n F_t(du_j)$$

$$\le K(\psi)e^{-2ax}(E_t[e^{2aX}])^n.$$

Fix $\eta \in (0, \frac{a}{L})$. Then for all $t \in I$,

$$\sum_{n\le \eta x}|J_n(t,x)| \le K(\psi)e^{-2ax}\sum_{n=1}^{\eta x} E_t[e^{2aX}]^n$$

$$\le \frac{K(\psi)e^{-2ax+L\eta x}}{1-e^{-L}} = \frac{K(\psi)e^{-ax}}{1-e^{-L}}.$$

On the other hand, by the selection of $I$, it is seen that there is $b > 0$, such that $l := \sup_{t\in I}\log E_t[e^{-bX}] < 0$. By Chernoff's inequality, for $n \ge 1$,

$$|J_n(t,x)| \le K(\psi)P_t\left(\sum_{j=1}^n X_j \le x\right) \le K(\psi)e^{bx}E_t[e^{-bX}]^n \le K(\psi)e^{bx+ln}.$$

By choosing $M \gg 0$, (4.13) is thus proved.

Let $A_x = [\eta x, Mx]$. We need to bound $\sum_{n\in A_x} J_n(s,t,x)$. For $n \ge \eta x$, let

$$J'_n(s,t,x) = \int \psi_{x-w_n}(t)\mathbf{1}\{w_k \le x, \forall k \in (\eta x, n]\}\prod_{j=1}^n E_t[e^{isY}|u_j]F_t(du_j).$$

With the same $a > 0$ and $L > 0$ as in the bound for $\sum_{n\le \eta x}|J_n|$, for $t \in I$,

$$|J_n(s,t,x) - J'_n(s,t,x)| \le \int \psi_{x-w_n}(x)\mathbf{1}\{w_k \ge x, \exists k \le \eta x\}\prod_{j=1}^n F_t(du_j)$$

$$(4.14) \qquad \le K(\psi)\sum_{k=1}^{\eta x} P_t\left(\sum_{j=1}^k X_j \ge x\right)$$

$$\leq K(\psi)e^{-2ax}\sum_{k=1}^{\eta x}e^{kL} = \frac{K(\psi)e^{-ax}}{1-e^{-L}}.$$

We next find an $a > 0$ such that $\sup_{t \in I}\sum_{n \in A_x} J'_n(s,t,x) = O(e^{-axs^2})$ and then use the above approximation to bound $\sup_{t \in I}\sum J_n(s,t,x)$. Note that

$$J'_n(s,t,x) = \int J_{n-\lfloor\eta x\rfloor}(s,t,x-w_{\lfloor\eta x\rfloor})\prod_{j=1}^{\eta x}E_t[e^{isY}|u_j]F_t(du_j).$$

Fix $k$ as in (1.17). Let $\xi = \sum_{j=1}^k Y_j$, $Z = \sum_{j=1}^k X_j$. In order to obtain the desired bound, we express $J'_n(s,t,x)$ in terms of $E_t[e^{is\xi}|Z]$. By (1.17), for each $t \in I$, $E_t|E_t[e^{is\xi}|Z]| < 1$, implying that $E_t[\sigma_Z(t)] > 0$, where $\sigma_z(t) = \text{Var}_t(\xi|Z=z)$. Note that if $\sigma_z(t) > 0$ for one $t \in I$, then by the smoothness of $P_t$, $v(z) := \inf_{t \in I}\sigma_z(t) > 0$. Therefore, $\Pr\{v(Z) > 0\} > 0$. On the other hand, by (1.18), $\sup_{t \in I} E[e^{t\xi+\eta|\xi|}|Z] < \infty$ a.s. Thus there exist $r > 0$ and $R > 1$ such that $p := \inf_{t \in I} P_t(Z \in \Gamma) > 0$, where $\Gamma = \{z : r \leq v(z), \sup_{t \in I} E_t[|\xi|^3|Z=z] \leq R\}$. Let $m = \lfloor\frac{\lfloor\eta x\rfloor}{k}\rfloor$, $Z_j = \sum_{s=1}^k X_{(j-1)k+s}$, $\xi_j = \sum_{s=1}^k Y_{(j-1)k+s}$. Then

$$\sum_{n \in A_x} J'_n(s,t,x) = \int \underbrace{\left(\sum_{k \leq (M-\eta)x} J_k(s,t,x-w_{\lfloor\eta x\rfloor})\right)}_{V}$$

$$\times \prod_{j=1}^m E_t[e^{is\xi_j}|Z_j = z_j]F_t^{k*}(dz_j)$$

$$\times \prod_{j=mk+1}^{\lfloor\eta x\rfloor} E_t[e^{isY_l}|X_j = u_j]F_t(du_j).$$

Since $|V| \leq \sum J_k(t, x-w_{\lfloor\eta x\rfloor}) = \phi_{x-w_{\lfloor\eta x\rfloor}}(t) \leq \sup_{t,x}\phi_x(t)$,

(4.15)
$$\left|\sum_{n \in A_x} J'_n(s,t,x)\right| \leq \sup_{t,x}\phi_x(t)\int \prod_{j=1}^m |E_t[e^{is\xi}|Z=z_j]|F_t^{k*}(dz_j)$$

$$\leq \sup_{t,x}\phi_x(t)\left[1-p+p\sup_{z \in \Gamma}|E_t[e^{is\xi}|Z=z]|\right]^m.$$

By Lemma 4.2, $\sup_{t,x}\phi_x(t) < \infty$. Let $\mu_t(z) = E_t[\xi|Z=z]$. Then by [2], Proposition 8.44, for $|s| \ll 1$, $z \in \Gamma$ and $t \in I$,

$$|E_t[e^{is\xi}|Z=z]| = |E_t[e^{is(\xi-\mu_t(z))}|Z=z]|$$

$$\leq 1 - \frac{rs^2}{2} + 4R|s|^3 \leq 1 - \frac{rs^2}{4}.$$

Because $\frac{m}{x} \to \frac{\eta}{k}$ as $x \to \infty$, it is seen that there exists $\alpha > 0$ such that for $|s| \ll 1$ and $x \gg 0$, $|\sum_{n \in A_x} J_n'(s,t,x)| \leq Ce^{-\alpha x s^2}$, with $C > 0$ an unspecified constant that may vary from appearance to appearance. Together with (4.14), this implies that $|\sum_{n \in A_x} J_n(s,t,x)| \leq C(xe^{-ax} + e^{-\alpha x s^2})$. Then by (4.13), there exists $\delta > 0$ such that for $|s| \leq \delta$, $t \in I$ and $x \gg 0$, $|G_x(t+is)| \leq C(e^{-\eta x} + xe^{-ax} + e^{-\alpha x s^2})$. This then completes the proof of (4.12). $\square$

PROOF OF (4.11). The formula will follow if $\sup_{t \in I, \eta < |s| < \lambda} |G_x(t+is)| = o(x^{-1/2})$, where $G_x(\cdot)$ is defined as in (4.12). From the proof of (4.10), $\sup_{t \in I, \eta < |s| < \lambda} |G_x(t+is) - \sum_{n \in A_x} J_n'(s,t,x)| = O(e^{-ax})$ for some $a > 0$. On the other hand, from (4.15), we obtain $|\sum_{n \in A_x} J_n'(s,t,x)| \leq Cf(s)^{\eta x/k}$, where $C$ is a constant and $f(s) = \sup_{t \in I} E_t|E_t[e^{is\xi}|Z]|$. As $E|E[e^{is\xi}|Z]| < 1$, it can be seen that $f(s) < 1$. The proof then follows. $\square$

## 5. Consequences of the uniform exact LDP.

PROOF OF COROLLARY 1.3. Let $\nu = h'(\tau_0)$. For $x \gg 0$, there exists a unique $t_x \in I$ with $h'(t_x) = \nu_x := \frac{x\nu + b}{x+a}$. As $t_x \to \tau_0$, by the uniform convergence,

$$\Pr\{W(x+a) \geq xh'(\tau_0) + b\} \sim \frac{\phi(t_x)}{t_x\sqrt{2\pi x h''(t_x)}} e^{-(x+a)h^*(\nu_x)}.$$

Since $t_x \to \tau_0 > 0$, $\frac{\phi(t_x)}{t_x\sqrt{2\pi x h''(t_x)}} = (1+o(1))\frac{\phi(\tau_0)}{\tau_0\sqrt{2\pi x h''(\tau_0)}}$. On the other hand,

$$-(x+a)h^*(\nu_x) + xh^*(\nu)$$
$$= x[h^*(\nu) - h^*(\nu_x)] - ah^*(\nu_x)$$
$$= h^{*\prime}(z)\frac{x(a\nu - b)}{x+a} - ah^*(\nu_x) \to h^{*\prime}(\nu)(a\nu - b) - ah^*(\nu),$$

with $z = z(x)$ between $\nu$ and $\nu_x$. Since $h^{*\prime}(\nu) = \tau_0$ and $h^*(\nu) = \nu\tau_0 - h(\tau_0)$, the limit is $ah(\tau_0) - b\tau_0$. This then completes the proof. $\square$

PROOF OF COROLLARY 1.4. Use the same $g_x(t)$, $k_x(t)$ and $F_t(dx)$ as in the proof of Theorem 1.2. On the other hand, define

$$\bar{M}_x(t) := E[e^{t\bar{W}(x)}], \qquad \bar{\phi}_x(t) = \frac{\bar{M}_x(t)}{e^{h(t)x}},$$

$$\bar{\psi}_x(t) = \frac{1}{e^{h(t)x}} \int_x^\infty g_u(t) F(du).$$

As in the proof of Theorem 1.2, let $\bar{M}_n(t) = M_{x_n}(t)$, $\bar{\Lambda}_n(t) = \bar{\Lambda}_{x_n}(t)$ and

$$\Pr\{\bar{A}_n^* \in du - x_n \nu_n\} = \frac{\exp\{t_n^* u\}}{\bar{M}_n(t_n^*)} \Pr\{\bar{W}(x_n) \in du\},$$

$$\bar{U}_n = \frac{\bar{A}_n^*}{\sqrt{x_n \bar{\Lambda}_n''(t_n^*)}}.$$

Clearly, Lemma 4.1 still applies. The proof then follows from Lemma 5.1 below and almost the same steps as in the proof of Theorem 1.2. □

LEMMA 5.1. *Under the conditions of Corollary 1.4, Lemma 4.2 still holds if $\Lambda_x(t)$, $\phi_x(t)$ and $\phi(t)$ are replaced by $\bar{\Lambda}_x(t) = \frac{1}{x} \log \bar{M}_x(t)$, $\bar{\phi}_x(t)$ and $\bar{\phi}(x)$, respectively, and equations (4.10) and (4.11) still hold if $M_n(t)$, $\Lambda_n(t)$ and $U_n$ are replaced by $\bar{M}_n(t)$, $\bar{\Lambda}_n(t)$ and $\bar{U}_n$, respectively.*

PROOF OF PROPOSITION 1.1. Let the span of $Y - f(X)$ be $d$. Since $X$ has a continuous density, each connected component of sppt$(X)$ has a nonempty interior. We first show that if $f(x)$ is continuous but not affine on a nonempty $(a,b) \subset \text{sppt}(X)$, then $(X,Y)$ satisfies condition (1.17) with $k = 2$. Let $(\xi_1, \zeta_1)$ and $(\xi_2, \zeta_2)$ be independent, such that $\Pr\{\xi_i \in dx\} = \Pr\{X \in dx | X \in [a,b]\}$ and $\Pr\{\zeta_i \in dy | \xi_i = x\} = \Pr\{Y \in dy | X = x\}$, $\forall x \in [a,b]$. It suffices to show that

(5.1) $$E|E[e^{is(\zeta_1 + \zeta_2)}|\xi_1 + \xi_2]| < 1 \qquad \forall s \neq 0.$$

Assume that (5.1) is not true for $s \neq 0$. Then for a.e. $u \in (a,b)$, there exists $h$ such that $\Pr\{\zeta_1 + \zeta_2 \in h + \frac{1}{2\pi s}\mathbb{Z} | \xi_1 + \xi_2 = 2u\} = 1$. Since $\zeta_1 + \zeta_2 \in f(\xi_1) + f(\xi_2) + d\mathbb{Z}$ a.s., $\Pr\{(h + \frac{1}{2\pi s}\mathbb{Z}) \cap (f(\xi_1) + f(\xi_2) + d\mathbb{Z}) \neq \varnothing | \xi_1 + \xi_2 = 2u\} = 1$, implying that for a.e. $x \in (2u-b, 2u-a) \cap (a,b)$, $f(2u-x) + f(x) \in h + \frac{1}{2\pi s}\mathbb{Z} + d\mathbb{Z}$. Since $f(2u-x) + f(x)$ is continuous in $x$ and $\frac{1}{2\pi s}\mathbb{Z} + d\mathbb{Z}$ is discrete, there exists $c$ such that $f(2u-x) + f(x) \equiv c$. In particular, letting $x = u$ yields $c = 2f(u)$. It follows that for $x, y \in (a,b)$, $f(x) + f(y) = 2f(\frac{x+y}{2})$ and hence $f$ is affine on $(a,b)$, a contradiction.

It remains to show (5.1) for the case where $f$ is piecewise affine. By the assumption, there exist disjoint intervals $(a_1, b_1), (a_2, b_2) \subset \text{sppt}(X)$ and constants $c_1 \neq c_2$, $d_1$ and $d_2$ such that $f(x) = c_i x + d_i$ on $(a_i, b_i)$. Let $(\xi_1, \zeta_1)$ and $(\xi_2, \zeta_2)$ be independent such that $\Pr\{\xi_i \in dx\} = \Pr\{X \in dx | X \in I_i\}$ and $\Pr\{\zeta_i \in dy | \xi_i = x\} = \Pr\{Y \in dy | X = x\}$. Now sppt$(\xi_1 + \xi_2) = J = [a_1 + a_2, b_1 + b_2]$. Given $\xi_1 + \xi_2 = u \in J^o$, $Z := f(\xi_1) + f(\xi_2) \stackrel{d}{=} (c_2 - c_1)\eta + c_1 u + d_1 + d_2$, where $\Pr\{\eta \in dx\} = \Pr\{\xi_2 \in dx | \xi_2 \in I\} = \Pr\{X \in dx | X \in I\}$, with $I = (u - b_2, u - a_1) \cap (a_2, b_2) \neq \varnothing$. Thus $Z$ has a density, yielding $|E[e^{is(\zeta_1 + \zeta_2)}|\xi_1 + \xi_2 = u]| = |E[e^{isZ}|\xi_1 + \xi_2 = u]| < 1$ for all $s \neq 0$. □

## 6. Auxiliary technical details.

6.1. *Proofs of auxiliary results for the renewal measures.*

PROOF OF LEMMA 3.1. By condition (1.4), $\hat{p}(z) = E_p[e^{izX}]$ are analytic and equicontinuous in $D_\tau$. Fix $r_0 \in (0, \tau)$. Let $\bar{\mu}_p = E_p|X|$ and $\nu_p = \frac{1}{2}E_p[X^2]$. For $z$ with $|z| \leq r_0$,

$$
\begin{aligned}
|\hat{p}(z) - 1 - \mu_p z - \nu_p z^2| &\leq \sum_{n=3}^{\infty} \frac{|z|^n E_p|X|^n}{n!} \\
&= E[e^{|zX|}] - 1 - \bar{\mu}_p|z| - \nu_p|z|^2 \\
&\leq E_p[|X|^3 e^{r_0|X|}]|z|^3 \leq C M_\tau |z|^3,
\end{aligned}
\tag{6.1}
$$

where $C \geq \sup_{x \geq 0} x^3 e^{(r_0 - \tau)x}$. By (1.4), $\sup_p \nu_p < \infty$. Combined with (1.3), this implies that there exists $\delta > 0$ such that $|\hat{p}(z) - 1| \geq \delta|z|$ if $|\Re(z)| \leq \delta$ and $|\Im(z)| \leq \delta$.

Assuming the first statement in Lemma 3.1 were false, there would be $t_n \to 0$, $\theta_n \in \mathbb{R}$ and $p_n \in \mathcal{M}$, such that $\theta_n - it_n \neq 0$ and $\hat{p}_n(\theta_n - it_n) = 1$. From the previous paragraph, $|\theta_n| \geq \delta$. Then by $|\hat{p}_n(\theta_n - it_n) - \hat{p}_n(\theta_n)| \leq E_{p_n}[e^{|t_n X|}] - 1 \leq M_\tau^{|t_n|/\tau} - 1 \to 0$, $\hat{p}_n(\theta_n) \to 1$. On the other hand, by (1.6) and (3.7), $\hat{p}_n(\theta_n) = (1 - \lambda_p)\hat{\Phi}_p(\theta_n) + \lambda_p \hat{\Psi}_p(\theta_n) \leq C\theta_n^{-2} + \gamma$. Therefore, $\liminf C\theta_n^{-2} \geq 1 - \gamma > 0$ and the $\theta_n$ are bounded. By (1.4), $\mathcal{M}$ is tight. Then there exist a probability measure $p_0$ and an $s \neq 0$ such that, say, $p_n \xrightarrow{w} p_0$ and $\theta_n \to s$. Note that $|s| \geq \delta$. Because $\hat{p}$ are equicontinuous on $D_\tau$, by $\hat{p}_n(\theta_n - it_n) = 1$, we obtain $\hat{p}_0(s) = 1$. Therefore, $p_0$ is concentrated on a lattice $L$ of span $|s|$.

Given $\varepsilon > 0$ and $\eta > 0$, let $L_\varepsilon = \bigcup_{x \in L}(x - \varepsilon, x + \varepsilon)$. There exists $R > T$ such that $p_n((-R, R) \cap L_\varepsilon) \geq 1 - \eta$, $\forall n \gg 1$. Then by (1.5), (1.6) and $\mathrm{sppt}(\Phi_p) \subset (-T, T)$,

$$
\gamma + \sup_{p \in \mathcal{M}} \Phi_p((-T, T) \cap L_\varepsilon) \geq 1 - \eta.
\tag{6.2}
$$

Because $\phi_p \in C_0^2((-T, T))$, from $\sup_{x,y}|\phi_p^{(i)}(x) - \phi_p^{(i)}(y)| \leq \int_{-T}^T |\phi_p^{(i+1)}|$, $i = 0, 1$, we obtain $A := \sup_{p \in \mathcal{M}} \|\phi_p\|_\infty < \infty$ and by (6.2), $\gamma + A\ell((-T, T) \cap L_\varepsilon) \geq 1 - \eta$. Letting $\varepsilon \to 0$ yields $\gamma \geq 1 - \eta$. Since $\eta$ is arbitrary, we arrive at a contradiction to (1.6). This proves the first statement in Lemma 3.1.

To show that $L_\eta < \infty$ for some $\eta \in (0, \tau)$, given $p \in \mathcal{M}$,

$$
|K_p(z)| = \left|\frac{1}{1 - \hat{p}(z)} - \frac{1}{-i\mu_p z}\right| = \frac{|\hat{p}(z) - 1 - i\mu_p z|}{\mu_p|\hat{p}(z) - 1||z|}.
$$

By (6.1), $|\hat{p}(z) - 1 - i\mu_p z| \leq \nu_p|z|^2 + CM_\tau|z|^3$ and $|\hat{p}(z) - 1| \geq \mu_p|z| - \nu_p|z|^2 + CM_\tau|z|^3$. As $\inf_p \mu_p > 0$ and $\sup_p \nu_p < \infty$, there exists $\delta > 0$ with

$\sup_{p \in \mathcal{M}, |\Re z|, |\Im z| \leq \delta} |K_p(z)| < \infty$. Thus if $L_\eta = \infty$ for any $\eta \in (0, \tau)$, then there exist $t_n \to 0$, $\theta_n \in \mathbb{R}$ and $p_n \in \mathcal{M}$ such that $|\theta_n| \geq \delta$ and $\hat{p}_n(\theta_n - it_n) - 1 \to 0$. The same argument for the first statement then applies to yield a contradiction. This completes the proof of Lemma 3.1. $\square$

6.2. *The moment generating function of the random sum.* This part provides technical results on $M_x(t)$ for the proof of Lemma 4.2. Recall that $\beta_X = \limsup_x \frac{1}{x} \log \Pr\{X > x\}$.

LEMMA 6.1. *The following statements* (1)–(4) *are true.* (1) *For any $t$ with $\mathcal{D}_t \neq \varnothing$, $E[e^{tY-sX}]$ is strictly convex in $s \in \mathcal{D}_t$. Furthermore, given $t_1, t_2$ and $\lambda \in (0,1)$, we have $(1-\lambda)\mathcal{D}_{t_1} + \lambda \mathcal{D}_{t_2} \subset \mathcal{D}_{(1-\lambda)t_1 + \lambda t_2}$.* (2) *For any $t \in \mathbb{R}$, $\mathcal{D}_t$ is a closed interval and $h(t) > -\infty$. In particular, if $h(t) < \infty$, then $h(t) \in \mathcal{D}_t$.* (3) *$h$ is convex in $\mathcal{A} = \{t \in \mathbb{R} : h(t) < \infty\}$ and $h(0) = 0$.* (4) *$\mathcal{A}$ is an interval with $0 \in \mathcal{A}$. In particular, for any $t \in \mathcal{A}$, if $\mathcal{D}_t^o \neq \varnothing$, then $\mathcal{D}_s^o \neq \varnothing$, $\forall s \in \mathcal{A}^o$.*

PROOF. (1) follows from the strict convexity of $e^x$, the assumption that $X$ is nondegenerate and Hölder's inequality. (2) By (1) and Fatou's lemma, $\mathcal{D}_t$ is a closed interval. If $h(t) = -\infty$, then for any $M \gg 0$, $E[e^{tY+MX}] \leq 1$, implying that $\text{ess sup } X \leq 0$, which contradicts (1.1). (3) By (1), $h$ is convex. Clearly, $h(0) \leq 0$. If $h(0) < 0$, then $E[e^{|h(0)|X}] \leq 1$ and by the strict convexity of $E[e^{sX}]$, there exists $a > 0$ such that $E[e^{aX}] < 1$. For $x \gg 0$, by Chernoff's inequality, $\sum \Pr\{\sum_{i=1}^n X_i > x\} \leq e^{-ax} \sum (E[e^{aX}])^n < 1$, implying that $\Pr\{N(x) = \infty\} > 0$, which contradicts (1.1). (4) follows easily from (1). $\square$

LEMMA 6.2 [Upper bound for exponential rate of $M_x(t)$]. *Let $t \in \mathbb{R}$ be such that $\beta_X < h(t) < \infty$ and $\mathcal{D}_t^o \neq \varnothing$. Then $M_x(t)$ is rightcontinuous in $x \geq 0$ and for $\varepsilon > 0$ with $h(t) + \varepsilon \in \mathcal{D}_t^o$, there exists $K_\varepsilon > 0$ such that*

$$(6.3) \qquad M_x(t) \leq K_\varepsilon e^{(h(t)+\varepsilon)x} \qquad \forall x > 0.$$

*Furthermore, if either $h(t) \geq 0$ or $\varepsilon < |h(t)|$, then one can set $K_\varepsilon = \frac{H(h(t))}{1-\rho(t,\varepsilon)}$, where $H(a) = E[e^{-(0 \wedge a)X}] + 1$ and $\rho(t, \varepsilon) = E[e^{tY-(h(t)+\varepsilon)X}]$.*

PROOF. For brevity, define $a_1^n = \sum_{i=1}^n a_i$. First, suppose that $h(t) \geq 0$. Fix $\varepsilon > 0$ with $h(t) + \varepsilon \in \mathcal{D}_t^o$. By Lemma 6.1, $E[e^{tY-(h(t)+\varepsilon)X}] < 1$. Then

$$(6.4) \quad \begin{aligned} e^{tW(x)} &= \sum_{k=0}^\infty e^{tY_1^k} \mathbf{1}\{N(x) = k\} \leq \sum_{k=0}^\infty e^{tY_1^k} \mathbf{1}\{X_1^k \leq x\} \\ &\leq \sum_{k=0}^\infty e^{tY_1^k} e^{(h(t)+\varepsilon)(x-X_1^k)} = e^{(h(t)+\varepsilon)x} \sum_{k=0}^\infty \xi_k, \end{aligned}$$

where $\xi_k = \prod_{j=1}^k \exp\{tY_j - (h(t)+\varepsilon)X_j\}$. Let $\zeta_x = e^{(h(t)+\varepsilon)x}\sum \xi_k$. Then

(6.5) $$E[\zeta_x] = e^{(h(t)+\varepsilon)x} \sum_{k=0}^{\infty} \rho(t,\varepsilon)^k < \infty.$$

Equations (6.4) and (6.5) show that $M_x(t) \leq K_\varepsilon e^{(h(t)+\varepsilon)x}$ and, given $x_0 > 0$, $e^{tW(x)} \leq \zeta_{x_0}$, $\forall x \in (0, x_0)$. When $y \downarrow x$, $N(y) \downarrow N(x)$ and $W(y) \to W(x)$ a.s. Then by dominated convergence, $M_x(t)$ is rightcontinuous in $x$.

If $\beta_X < h(t) < 0$, then we show that for $0 \leq \varepsilon < -h(t)$, $E[e^{-(h(t)+\varepsilon)X}] < \infty$. Fix $\delta \in (0, h(t) - \beta_X)$. For $x \gg 1$, $\Pr\{X > x\} \leq e^{-(|h(t)|+\delta)x}$. Then

$$\int_0^\infty e^{-(h(t)+\varepsilon)x} F(dx) \leq \int_0^\infty e^{-h(t)x} F(dx)$$
$$\leq 1 + |h(t)| \int_0^\infty e^{|h(t)|x} \Pr\{X > x\} dx < \infty$$

and hence $E[e^{-(h(t)+\varepsilon)X}\mathbf{1}\{X \geq 0\}] < \infty$. On the other hand, because $h(t)+\varepsilon < 0$, $E[e^{-(h(t)+\varepsilon)X}\mathbf{1}\{X < 0\}] < 1$. Therefore, $E[e^{-(h(t)+\varepsilon)X}] \leq E[e^{-h(t)X}] + 1 < \infty$. Fix $\varepsilon \in (0, -h(t))$ with $h(t) + \varepsilon \in \mathcal{D}_t^o$. Then

$$e^{tW(x)} = \sum_{k=0}^{\infty} e^{tY_1^k}\mathbf{1}\{N(x) = k\} \leq \sum_{k=0}^{\infty} e^{tY_1^k}\mathbf{1}\{X_1^{k+1} > x\}$$

(6.6) $$\leq \sum_{k=0}^{\infty} e^{tY_1^k} e^{-(h(t)+\varepsilon)(X_1^{k+1}-x)}$$

$$= e^{(h(t)+\varepsilon)x} \sum_{k=0}^{\infty} \xi_k e^{-(h(t)+\varepsilon)X_{k+1}}.$$

Note that $\xi_n$ is independent of $X_{n+1}$. An argument similar to that for $h(t) \geq 0$ can then be applied. The expression for $K_\varepsilon$ follows from (6.5) and (6.6). □

LEMMA 6.3 [Upper bound for exponential rate of $M_x^{(n)}(t)$]. *Suppose that there exists $I = (a,b)$ such that $\beta_X < h(t) < \infty$ on $I$. Also, suppose that there exist $\tau \in I$ and $\eta > 0$ such that $\mathcal{D}_\tau^o \neq \varnothing$ and $E[e^{tY-h(t)X+\eta(|X|+|Y|)}] < \infty$ for $t \in I$. Then for any $x \geq 0$, $M_x(t) \in C^\infty(I)$ and there exists $K = K(n,t,\varepsilon)$, such that $|\frac{\partial^n M_x(t)}{\partial t^n}| \leq Ke^{(h(t)+\varepsilon)x}$, $\forall t \in I$, $\varepsilon > 0$. Furthermore, if $h(t)+\varepsilon \in \mathcal{D}_t^o$ and either $h(t) \geq 0$ or $\varepsilon < |h(t)|$, then one can set $K = H(h(t)) \times H_1(h(t)) \times P_n(\rho(t,\varepsilon)) \times (1 - \rho(t,\varepsilon))^{-n}$, where $H(\cdot)$ is given in Lemma 6.2, $H_1(a) = CE[e^{tY-h(t)X+\eta(|X|+|Y|)}]$ with $C$ an absolute constant and $P_n$ a polynomial of degree $n$ with absolute constant coefficients. Finally, $\frac{\partial^n M}{\partial t^n} \in C(I \times \mathbb{R})$.*

PROOF. Because $\mathcal{D}_\tau^o \neq \varnothing$, by Lemma 6.1, $\mathcal{D}_t^o \neq \varnothing$, $\forall t \in I$. Then by Lemma 6.2, $\forall x \geq 0$, $M_x(t) < \infty$ on $I$, implying that $M_x(t) \in C^\infty(I)$ with

$M_x^{(n)}(t) = E[W(x)^n e^{tW(x)}]$. For $n \geq 1$, by $|\sum_{i=1}^k x_i|^n \leq k^{n-1} \sum_{i=1}^k |x_i|^n$, if $h(t) \geq 0$ and $0 < \varepsilon \ll 1$, then, as in (6.4),

$$|W(x)|^n e^{tW(x)} \leq e^{(h(t)+\varepsilon)x} \sum_{k=0}^\infty k^{n-1} \left(\sum_{i=1}^k |Y_i|^n\right) \xi_k.$$

Let $C = \sup_{x \geq 0} x^n e^{-\eta x}$ and $\zeta_i = C e^{tY_i - h(t)X_i + \eta(|X_i|+|Y_i|)}$. Since $\zeta_i$ are i.i.d. with $E\zeta_1 < \infty$ and $|Y_i|^n e^{tY_i - (h(t)+\varepsilon)X_i} \leq \zeta_i$, we have $E[|W(x)|^n e^{tW(x)}] \leq Ke^{(h(t)+\varepsilon)x}$, where $K = E\zeta_1 \sum_{k=1}^\infty k^n \rho(t,\varepsilon)^{k-1}$. The case $\beta_X < h(t) < 0$ is proved likewise. The rest of the proof follows that of Lemma 6.2. □

LEMMA 6.4 [Exact exponential rate of $M_x(t)$]. *Suppose that $\beta_X < h(t) < \infty$ and $E[e^{tY - h(t)X}] = 1$. Suppose that either (1) $\Pr\{X \geq 0\} = 1$, or (2) $\mathcal{D}_t^o \neq \varnothing$ and $E[e^{q(tY-(0\wedge h(t))X)}|X > 0] < \infty$, $\forall q > 0$. Then $\frac{1}{x} \log M_x(t) \to h(t)$.*

PROOF. By Lemma 6.2, we have $\limsup \frac{1}{x} \log M_x(t) \leq h(t)$. It remains to demonstrate that $\liminf \frac{1}{x} \log M_x(t) \geq h(t)$. First suppose that (1) is true. Then $\mathcal{D}_t = [h(t), \infty) \neq \varnothing$. Assume that the lower bound does not hold. Then there exist $\varepsilon > 0$ and $z_n \to \infty$ such that $f(z_n) := \frac{M_{z_n}(t)}{e^{(h(t)-\varepsilon)z_n}} \to 0$. Fix $\eta_n \downarrow 0$. Then $x_n := \sup\{x : \forall u \in [0,x), f(u) \geq \eta_n\} < \infty$ and

$$f(x_n) \geq \int_0^x \frac{g_u(t)}{e^{(h(t)-\varepsilon)u}} f(x_n - u) F(du) \geq \eta_n E[e^Z \mathbf{1}\{X \in (0, x_n]\}],$$

where $Z = tY - (h(t)-\varepsilon)X$. Since $M_x(t)$ is rightcontinuous, $f(x_n) \leq \eta_n$ and hence $E[e^Z \mathbf{1}\{X \in [0, x_n]\}] \leq 1$. Note that $x_n$ is increasing. If $x_n \uparrow \infty$, then $E[e^Z] \leq 1$, contradicting the definition of $h(t)$. If $x_n \to z < \infty$, then $M_{x_n}(t) = f(x_n) e^{(h(t)-\varepsilon)x_n} \to 0$, that is, $E[e^{tW(x_n)}] \to 0$, which implies that $\Pr\{|W(x_n)| > M\} \to 1$ for any $M > 0$. However, since $N(z) < \infty$ and $|Y| < \infty$ a.s., this is impossible. Thus the lower bound holds.

Suppose that (2) is true. Let $U_n = \exp\{tY_n - h(t)X_n\}$, $Z_n = U_1 \cdots U_n$ and $T(x) = N(x) + 1$. Fix $\varepsilon > 0$ with $h(t) + \varepsilon \in \mathcal{D}_t^o$. Define $\xi_n$ as in (6.4). Then

$$Z_n \mathbf{1}\{T(x) = n\} \leq Z_n \mathbf{1}\{T(x) \geq n\}$$
$$\leq \prod_{k=1}^n e^{tY_k - h(t)X_k} \times e^{\varepsilon(x - \sum_{i=1}^{n-1} X_i)} = e^{\varepsilon x} \xi_{n-1} U_n.$$

Since $E[\xi_{n-1}] = \rho(t,\varepsilon)^{n-1}$ and $E[U_n] = 1$, with $\rho(t,\varepsilon) < 1$ defined as in Lemma 6.2, there exists $M > 1$ such that for $x \gg 1$, $E[Z_{T(x)} \mathbf{1}\{T(x) \geq Mx\}] < \frac{1}{2}$. Let $X^* = X^*(x) = \sum_{i=1}^{T(x)} X_i - x$, $\zeta = \zeta(x) = \frac{1}{T(x)} e^{tY_{T(x)} - h(t)X^*}$. Then

$$E[e^{tW(x)}\zeta] \geq \frac{1}{Mx} E[e^{tW(x)} e^{tY_{T(x)} - h(t)X^*} \mathbf{1}\{T(x) < Mx\}] \geq \frac{e^{h(t)x}}{2Mx}.$$

Fix $p, q > 1$ with $\frac{1}{p} + \frac{1}{q} = 1$ and $pt \in \mathcal{A}^o$. Then $E[e^{ptW(x)}]^{1/p} E[\zeta^q]^{1/q} \geq \frac{e^{h(t)x}}{2Mx}$. Note that $M_x(pt) = E[e^{ptW(x)}]$. Also, note that $0 < X^* \leq X_n$ if $T(x) = n$. So by the assumption, $E[\zeta^q] \leq \sum n^{-q} E[e^{q(tY_n - (h(t) \wedge 0)X_n)} \mathbf{1}\{X_n > 0\}] < \infty$. Thus $\liminf \frac{1}{px} \log M_x(pt) \geq h(t)$. Replacing $pt$ by $t$, we get $\liminf \frac{1}{x} \log M_x(t) \geq ph(t/p)$, $\forall t \in \mathcal{A}^o$. Let $p \uparrow 1$. Since $h(t)$ is convex and hence continuous in $\mathcal{A}^o$, the lower bound follows. □

6.3. *Proofs of auxiliary results for the uniform exact LDP.*

PROOF OF LEMMA 4.1.  Recall the definition of $I_0$ and $\eta > 0$ given just before (4.2). By (1.18), there exists $a > 0$ such that

$$\sup_{t \in I_0} E_t[e^{a|X|}] \leq E[e^{\tau_0 Y - h(\tau_0) X + \eta_0(|X| + |Y|)}] < \infty.$$

Let $R > 0$ and $\pi(x) \in C_0^\infty(\text{sppt}(\Phi) \cap (-R, R))$ such that $0 \leq \pi(x) \leq 1$ and $\int \pi \phi > 0$. Now each $F_t$ has a subcomponent with density $k_x(t)\pi(x)\phi(x) \in C_0^2((-R, R))$. It is then seen that the conditions in Theorem 1.1 are satisfied. □

PROOF OF LEMMA 4.2.  By (4.2), $h \in C^\infty(I_0)$ and $E[e^{tY - h(t)X}] = 1$, one can choose $\varepsilon \in (0, \frac{\varepsilon_0 \wedge \eta}{4})$ and $I = (\tau_0 - \varepsilon, \tau_0 + \varepsilon)$ such that (1) $h(t) + \varepsilon \in \mathcal{D}_t^o$ and $|h(t) - h(\tau_0)| \leq \frac{\eta}{4}$ for every $t \in I$; (2) $\sup_{t \in I} \rho(t, a) < 1$, $\forall a \in (0, \varepsilon]$, where $\rho(t, a) = E[e^{tY - (h(t) + \varepsilon)X}]$; (3) $\inf_{t \in I} |h^{(n)}(t)| < \infty$, $n \geq 0$. By Lemma 6.3, $\forall x > 0$, $M_x(t) \in C^\infty(I)$ and $\forall 0 < \eta \ll 1$, there exists $K(\eta) > 0$ such that

(6.7) $\quad \sup_{t \in I} \{e^{-h(t)x}(|M_x(t)| + |M'_x(t)| + |M''_x(t)|)\} \leq K(\eta) e^{\eta x}$.

We first show (4.7). From the definitions, for $x \geq 0$,

(6.8) $\quad M_x(t) = \Pr\{X > x\} + \int_{-\infty}^x g_u(t) M_{x-u}(t) F(du),$

(6.9)
$$\phi_x(t) = \psi_x(t) + \int_{-\infty}^x k_u(t) \phi_{x-u}(t) F(du)$$
$$= \psi_x(t) + \int_{-\infty}^x \phi_{x-u}(t) F_t(du).$$

From (6.7), $\sup_{x \geq 0, t \in I}\{e^{-\varepsilon x} \phi_x(t)\} < \infty$, $\forall \varepsilon \ll 1$. On the other hand, by the selection of $I_0$, $\beta_X < \inf_{I_0} h \leq \inf_I h$. Since $\psi_x(t) = e^{-h(t)x} \Pr\{X > x\}$, there exists $\varepsilon > 0$ such that $\sup_{t \in I} \psi_x(t) = o(e^{-\varepsilon x})$. We now apply Lemma 4.1 and Corollary 1.2 to (6.9). Then

$$\phi_x(t) \to \frac{1}{E_t[X]} \left[ \int_0^\infty \psi_x(t) \, dx - \int_{-\infty}^0 \left( \int_0^{|x|} \phi_u(t) \, du \right) F_t(dx) \right]$$

uniformly exponentially fast for $t \in I$. The right-hand side is $\phi(t)$, so (4.7) holds.

It remains to show (4.6). Because $\Lambda_x(t) = \frac{1}{x} \log M_x(t) = \frac{1}{x} \log \phi_x(t) + h(t)$,

$$\Lambda'_x(t) = \frac{\phi'_x(t)}{x\phi_x(t)} + h'(t), \qquad \Lambda''_x(t) = \frac{\phi''_x(t)\phi_x(t) - \phi'_x(t)^2}{x\phi_x(t)^2} + h''(t).$$

Therefore, in order to obtain (4.6), it suffices to show that

(6.10) $$\limsup_{x \to \infty} \sup_{t \in I} \left| \frac{\phi_x^{(n)}(t)}{\phi_x(t)} \right| < \infty, \qquad n = 1, 2.$$

To show (6.10) for $\phi'_x(t)$, we first show that for $x > 0$ and $t \in I$,

(6.11) $$\phi'_x(t) = \psi'_x(t) + \int_{-\infty}^{x} (k_u(t)\phi_{x-u}(t))' F(du),$$

that is, $\frac{d}{dt} \int_{-\infty}^{x} k_u(t)\phi_{x-u}(t) F(du) = \int_{-\infty}^{x} (k_u(t)\phi_{x-u}(t))' F(du)$. It suffices to verify that for $0 < \varepsilon \ll 1$, $\int_{-\infty}^{x} \sup_{|s-t|\le\varepsilon} |(k_u(s)\phi_{x-u}(s))''| F(du) < \infty$. Indeed, as $k_x(t) = e^{-h(t)x} E[e^{tY}|X=x]$ and $\phi_x(t) = e^{-h(t)x} M_x(t)$, by (6.7), it can be seen that for $\varepsilon \ll 1$, there exists $C$ such that for $n \le 2$, $t \in I$ and $x \ge 0$, $|k_x^{(n)}(t)| \le Ce^{-h(t)x} E[e^{tY+\varepsilon|Y|}|X=x]$ and $\sup_I |\phi_x^{(n)}(t)| \le Ce^{\varepsilon x}$. By (1.18), fix $\eta > 0$ and $\varepsilon \in (0, \frac{\eta}{2})$ such that $Z := e^{tY-h(t)X+\eta(|X|+|Y|)} \in L^1$ and $\sup_{|s-t|\le\varepsilon} |h(s) - h(t)| < \frac{\eta}{2}$. It follows that $\sup_{|s-t|\le\varepsilon} |(k_u(s)\phi_{x-u}(s))''| \le 4C^2 e^{\varepsilon x} [Z|X=u] \in L^1(F)$ and hence (6.11) holds.

Now let $w_x(t) = \int_{-\infty}^{x} k'_u(t)\phi_{x-u}(t) F(du)$. By (6.11), $\phi'_x(t) = z_x(t) + \int_{-\infty}^{x} \phi'_{x-u}(t) F(du)$, where $z_x(t) = \psi'_x(t) + w_x(t)$. As $\int_{-\infty}^{\infty} k_u(t) F(du) \equiv 1$, by an argument similar to that above, $\int_{-\infty}^{\infty} k'_u(t) F(du) = 0$. Then

$$|w_x(t)| \le \int_{-\infty}^{x} |k'_u(t)| |\phi_{x-u}(t) - \phi(t)| F(du) + \phi(t) \int_x^{\infty} |k'_u(t)| F(du).$$

From (4.7), $\sup_{t \in I} |\phi_x(t) - \phi(t)| \le Ce^{-ax}$ for $0 < a \ll 1$. Combining the above results, we have

$$|w_x(t)| \le C \int_{-\infty}^{x} e^{-h(t)u - a(x-u)} E[e^{tY+\eta|Y|}|X=u] F(du)$$

$$+ \phi(t) \int_x^{\infty} e^{-h(t)u+\eta(x-u)} E[e^{tY+\eta|Y|}|X=u] F(du)$$

$$\le Ce^{-(\eta \wedge \varepsilon)x}.$$

On the other hand, because $\beta_X < \inf_{t \in I} h(t)$, $\phi'_x(t) = -xh'(t)e^{-h(t)x} \times \Pr\{X > x\} \to 0$ uniformly exponentially fast on $I$. Thus $\sup_{t \in I} |z_x(t)| \to 0$ exponentially fast and Corollary 1.2 can be applied to (6.11) to obtain $\sup_I |\phi'_x(t) - \phi'(t)| \to 0$ exponentially fast. It also implies that $\phi \in C^1$ (cf. [20], Theorem 7.17).

The exponentially fast convergence $\sup_{t \in I} |\phi''_x(t) \to \phi''(t)| \to 0$ can be proven likewise. To complete the proof, we must show that, by shrinking $I$ if necessary, $\inf_{t \in I} |\phi(t)| > 0$. First, if $X \geq 0$, a.s., then this follows easily. Suppose that $\Pr\{X \geq 0\} < 1$. Since $\phi_x(\tau_0) > 0$, we have $\phi(\tau_0) = \lim_x \phi_x(\tau_0) \geq 0$. If $\phi(\tau_0) = 0$, then $\phi_x(\tau_0) \to 0$ exponentially fast in $x$, that is, there exists $a > 0$ such that $M_x(\tau_0)e^{-h(\tau_0)x} = o(e^{-ax})$. This implies that $\limsup_x \frac{1}{x} \log M_x(\tau_0) < h(\tau_0)$. On the other hand, by condition (1.19) and Lemma 6.4, $(1/x) \log M_x(\tau_0) \to h(\tau_0)$. The contradiction implies that $\phi(\tau_0) > 0$. Since $\phi$ is continuous around $\tau_0$, we can shrink $I$ to obtain $\inf_{t \in I} \phi(t) > 0$. This completes the proof of (6.10). □

PROOF OF LEMMA 4.3. Without loss of generality, let $d = 1$ and $\Pr\{Y > 0\} > 0$. Clearly, $W(x) \subset \mathbb{Z}$. By simple results from number theory, there exist $y_1, \ldots, y_s \in \text{sppt}(Y)$ such that for each large $m \in \mathbb{N}$, there exist $c_1, \ldots, c_s \in \mathbb{N}$ with $\sum c_i y_i = m$. Fix $\Gamma_i$, $i \leq s$, such that $\Pr\{X \in \Gamma_i\} > 0$ and $\Pr\{Y = y_i | X = x\} > 0$, $\forall x \in \Gamma_i$. Since the law of $X$ has a subcomponent with a density and $\Pr\{\sup_n \sum_{i=1}^n X_i = \infty\} = 1$, there exists $t_0$, such that the law of $\xi_0 := \sum_{i \leq t_0} X_i$ has a subcomponent with a density in $(0, \infty)$. Then there exist $k_0$ and $A \subset \text{sppt}(\xi_0)^o \cap (0, \infty)$, with $\ell(A) > 0$ and $\Pr\{\sum_{i \leq t_0} Y_i = k_0 | \xi_0 = x\} > 0$, $\forall x \in A$. By the property of measurable sets in $\mathbb{R}$, $\{x + y : x, y \in A\}$ contains an interval $(a, b) \neq \varnothing$. Let $k = 2k_0$, $t = 2t_0$. Then $(a, b) \subset \text{sppt}(\sum_{i \leq t} X_i)^o$ and $\Pr\{\sum_{i \leq t} Y_i = k | \sum_{i \leq t} X_i = x\} > 0$, $\forall x \in (a, b)$. Finally, there exists $I = [t, t+1] \cap \text{sppt}(X)$ with $t > 0$ such that $\Pr\{X \in I\} > 0$. Without loss of generality, assume that all $\Gamma_i$, $(a, b)$ and $I$ are in $[-R, R]$ for some $R > 0$.

Fix $m$, $p_1, \ldots, p_s$, $q_1, \ldots, q_s \in \mathbb{N}$ such that $\sum p_i y_i = m$, $\sum q_i y_i = m + 1$. Define $p = \sum p_i$, $q = \sum q_i$. We claim that for $x \gg 0$, there exists $n \in \mathbb{N}$ such that the events $E_1$ and $E_2$ defined as follows each has a positive probability. $E_1$ is the joint event of (1) $N(x) = nt + p$, (2) $X_{N(x)+1} \in I$ and (3) $\{Z_i := (X_i, Y_i), i \leq nt + p\}$ can be partitioned into $B_1, \ldots, B_s$, $C_1, \ldots, C_n$ such that $|B_j| = p_j$ with $X_i \in \Gamma_j$ and $Y_i = y_j$ for each $Z_i \in B_j$, $C_j = \{Z_{i_{j1}}, \ldots, Z_{i_{jt}}\}$, with $\sum_{l \leq t} X_{i_{jl}} \in (a, b)$ and $\sum_{l \leq t} Y_{i_{jl}} = k$. $E_2$ is defined likewise, with $q_i$ replacing of $p_i$. Event $E_1$ implies that $W(x) = nk + \sum p_i y_i = nk + m$, while $E_2$ implies that $W(x) = nk + m + 1$. Therefore, $nk + m, nk + m + 1 \in \text{sppt}(W(x))$ and $W(x)$ has span 1.

To verify the claim, let $T_0 = 1 + \max\{p, q\}$. For $x \gg 0$, there exists $n \in \mathbb{N}$ such that $T_0 R + na < x$ and $-T_0 R + nb > x$. This implies that for any $z_{i1}, \ldots, z_{ip_i} \in \Gamma_i$, $i \leq s$, and $z \in I$, there exists a nonempty open interval $J \in (a, b)$ such that the inequalities

$$\sum_{i,j} z_{ij} + \sum_{i=1}^n u_i \leq x, \qquad \sum_{i,j} z_{ij} + \sum_{i=1}^n u_i + z > x$$

are satisfied by any $u_1, \ldots, u_n \in J$. Therefore, the probability that all of the following events on $(X_1, Y_1), \ldots, (X_{nt}, Y_{nt})$ happen simultaneously is positive: (1) $\sum_{i,j} z_{ij} + \sum_i X_i < x$, (2) $\sum_{i,j} z_{ij} + \sum_i X_i + z > x$ and (3) $\forall i \leq n$, $\sum_{l \leq t} X_{it+l} \in (a,b)$ and $\sum_{l \leq t} Y_{it+l} = k$. Let $X_{nt+1}, \ldots, X_{nt+p}$ take values $x_{11}, \ldots, z_{1p_1}, \ldots, z_{s1}, \ldots, z_{sp_s}$, respectively, while let $X_{nt+1} = z$, and $Y_{nt+1}, \ldots, Y_{nt+p}$ take values $y_1, \ldots, y_s$, respectively, such that $Y_{nt+i} = y_j$ if $\sum_{h=1}^{j-1} p_h < j \leq \sum_{h=1}^{j} p_h$. Then $\sum_{i \leq nt+p} Y_i = nk + m$. Rearrange $(X_1, Y_1), \ldots, (X_{nt+p}, Y_{nt+p})$ so that those with negative values of $X_i$ appear first. This results in $N(x) = nt + p$ and $W(x) = nk + m$. Hence the claim holds for $E_1$. The case for $E_2$ can be treated similarly. $\square$

PROOF OF LEMMA 5.1. In order to show Lemma 4.2 with $\Lambda_x(t)$, $\phi_x(t)$ and $\phi(t)$ being replaced by $\bar{\Lambda}_x(t)$, $\bar{\phi}_x(t)$ and $\bar{\phi}(x)$, note that

$$\bar{\phi}_x(t) = \bar{\psi}_x(t) + \int_{-\infty}^x k_u(t) \bar{\phi}_{x-u}(t) F(du)$$
$$= \bar{\psi}_x(t) + \int_{-\infty}^x \bar{\phi}_{x-u}(t) F_t(du),$$

as opposed to (6.9). From the proof of Lemma 4.2, one can see that the lemma is implied by Lemma 4.1 and the fact that

(6.12) $\quad \sup_{t \in I} |\psi_x^{(n)}(t)| = o(e^{-ax}), \qquad n = 0, 1, 2,$ for some $a > 0$.

Indeed, (6.12) implies Lemmas 6.2–6.4, which are then combined with Lemma 4.1, (6.12) and Corollary 1.2 to obtain (4.6) and (4.7). For the current proof, Lemma 4.1 still holds. Therefore, we only need to check (6.12) with $\psi_x$ replaced by $\bar{\psi}_x$, as well as Lemmas 6.2–6.4. By Hölder's inequality,

$$|\bar{\psi}_x(t)| = \frac{1}{e^{h(t)x}} \int_x^\infty g(t,u) F(du) \leq \frac{1}{e^{h(t)x}} E[e^{qtY}]^{1/q} \Pr\{X > x\}^{1-1/q},$$

for all $q > 1$. Set $q \gg 0$ so that $(1 - \frac{1}{q})\beta_X < h(\tau_0) < \infty$. By the assumption of Corollary 1.4, $E[e^{qtY}]^{1/q} < \infty$. Therefore, $\sup_{t \in I} |\bar{\psi}_x(t)| = o(e^{-ax})$ for some $a > 0$. The exponential decay of $\sup_{t \in I} |\bar{\psi}_x^{(n)}(t)|$ for $n = 1, 2$ can be shown likewise. The proofs of Lemmas 6.2–6.4 for $\bar{M}_x(t)$ closely follow those for $M_x(t)$. We omit the detail for brevity.

From the proofs of (4.10) and (4.11), one can see that they are implied by conditions (1.17) and (1.18) as well as the fact that there exists $a > 0$ such that $\sup_{t \in I, x \geq 0} |e^{2ax} \psi_x(t)| < \infty$. Conditions (1.17) and (1.18) are still assumed here. It is not hard to show that there exists $a > 0$ such that $\sup_{t \in I, x \geq 0} |e^{2ax} \bar{\psi}_x(t)| < \infty$. By now repeating the proofs, it is seen that (4.10) and (4.11) still hold when $\bar{M}_n(t)$, $\bar{\Lambda}_n(t)$ and $\bar{U}_n$ replace $M_n(t)$, $\Lambda_n(t)$ and $U_n$. $\square$

## APPENDIX

PROOF OF (3.2) AND (3.3). Let $P$ be the probability measure that has density $p$. As in [2], $\forall r \in (0,1)$, let $q_r$ be the density of $\chi_p * \sum_{n=0}^{\infty} r^n P^{n*}$. Then $\hat{q}_r = \hat{\chi}_p \sum_{n=0}^{\infty} r^n \hat{p}^n = \frac{(1-\lambda)\hat{\bar{N}}_p \hat{\Phi}_p}{1-r\hat{p}}$. Since $(1-\lambda_p)\bar{N}_p$ is a probability measure, $(1-\lambda_p)|\hat{\bar{N}}_p| \leq 1$. By (3.7), $\hat{\Phi}_p \in L^1$. Therefore, $\hat{\chi}_p = (1-\lambda_p)\hat{\Phi}_p\hat{\bar{N}}_p \in L^1$. By $|1-r\hat{p}| \geq 1-r$, $\hat{q}_r \in L^1$. Given $M > 0$, by an inverse Fourier transform,

$$q_r(x) = \frac{1}{2\pi}\int_{-\infty}^{\infty} e^{-ix\theta}\frac{\hat{\chi}_p(\theta)}{1-r\hat{p}(\theta)}\,d\theta = \frac{1}{2\pi}\int_{-\infty}^{\infty}\Re\left(e^{-ix\theta}\frac{\hat{\chi}_p(\theta)}{1-r\hat{p}(\theta)}\right)d\theta$$

$$= \underbrace{\frac{1}{2\pi}\int_{-M}^{M}\Re\left(e^{-ix\theta}\frac{\hat{\chi}_p(\theta)}{1-r\hat{p}(\theta)}\right)d\theta}_{I_1} + \underbrace{\frac{1}{2\pi}\int_{|\theta|>M}\Re\left(e^{-ix\theta}\frac{\hat{\chi}_p(\theta)}{1-r\hat{p}(\theta)}\right)d\theta}_{I_2}.$$

Let $e^{-ix\theta}\chi_p(\theta) = f(\theta) + ig(\theta)$, $\hat{p}(\theta) = u(\theta) + iv(\theta)$. Then

$$I_1 = \int_{-M}^{M} f(\theta)\Re\left(\frac{1}{1-r\hat{p}(\theta)}\right)d\theta + \int_{-M}^{M}\frac{rg(\theta)v(\theta)}{|1-r\hat{p}(\theta)|^2}\,d\theta = J_1 + J_2.$$

For $\theta \neq 0$, $\hat{p}(\theta) \neq 1$. Therefore, $\lim_{r\uparrow 1} J_1 = \frac{\pi}{\mu_p} + \int_{-M}^{M} f(\theta)\Re(\frac{1}{1-\hat{p}(\theta)})\,d\theta$ (cf. [2], Lemma 10.11 or [24], page 330). By its definition, $\Psi_p$ has a finite mean, say $a$, and hence $\bar{N}_p$ has finite first moment $a\sum n\lambda_p^n$. Consequently, $\chi_p$ has a finite mean, yielding $g(\theta) \sim \mu_{\chi_p}\theta$ as $\theta \to 0$. On the other hand, $\hat{p}(\theta) = 1 + \mu_p\theta + o(\theta)$. By dominated convergence for $J_2$, as $r \uparrow 1$,

$$J_2 \to \int_{-M}^{M}\frac{g(\theta)v(\theta)}{|1-\hat{p}(\theta)|^2}\,d\theta \quad\Longrightarrow\quad I_1 \to \frac{\pi}{\mu_p} + \int_{-M}^{M}\Re\left(e^{-ix\theta}\frac{\hat{\chi}_p(\theta)}{1-\hat{p}(\theta)}\right)d\theta.$$

Because $\hat{p}(\theta) \to 0$ as $\theta \to \infty$, dominated convergence applies to $I_2$. This then completes the proof of (3.2). To now obtain (3.3), first let $x = 0$. Because $\chi_p$ has a finite first moment and $\hat{\chi}_p \in L^1$, it is seen that $\chi_p(\theta)/\theta \in L^1$. For $s > 0$, let $g(s) = \frac{1}{2\pi}\int_{-\infty}^{\infty}\frac{1}{\theta}(\int_{-\infty}^{\infty}\sin(t\theta)\chi_p(dt))e^{-\theta^2/2s^2}\,d\theta$. It is seen that $g(s)$ converges to the left-hand side of (3.3) as $s \to \infty$. On the other hand, by Fubini's theorem, $g(s) = \frac{1}{2\pi}\int_{-\infty}^{\infty} h(t,s)\chi_p(dt)$, where $h(t,s) = \int_{-\infty}^{\infty}\frac{1}{\theta}\sin(t\theta)e^{-\theta^2/2s^2}\,d\theta$. Observe that $\partial_t h = 2\pi\varphi_s(t)$, where $\varphi_s$ is the density of $N(0,\frac{1}{s})$. By $h(0,s) = 0$, $h(t,s) = 2\pi(\Pr\{Z \leq st\} - \frac{1}{2})$, where $Z \sim N(0,1)$. Therefore, $g(s) = \int_{-\infty}^{\infty}\Pr\{Z \leq st\}\chi_p(dt) - \frac{1}{2}$. Let $s \to \infty$ and apply dominated convergence to complete the proof. □

DEPARTMENT OF STATISTICS
UNIVERSITY OF CONNECTICUT
215 GLENBROOK ROAD, U-4120
STORRS, CONNECTICUT 06269
USA
E-MAIL: zchi@stat.uconn.edu